\documentclass[11pt]{amsart}

\usepackage[utf8]{inputenc} 
\usepackage{graphics,graphicx}
\usepackage{enumitem}
\usepackage{color}
\usepackage{verbatim}
\usepackage{float}
\usepackage{amsmath,amsthm, amscd, amssymb,amsfonts,mathtools}
\usepackage{wasysym}
\usepackage[all]{xy}
\usepackage{multirow}
\usepackage[active]{srcltx}
\usepackage{array}

\usepackage{listings}

\newcommand{\Dchaintwo}[3]{\xymatrix@C-4pt{\overset{#1}{\underset{\ }{\circ}}\ar
		@{-}[r]^{#2}
		& \overset{#3}{\underset{\ }{\circ}}}}

\makeindex

\newtheorem{theorem}{Theorem}[section]
\newtheorem{lemma}[theorem]{Lemma}
\newtheorem{corollary}[theorem]{Corollary}

\newtheorem{proposition}[theorem]{Proposition}
\newtheorem*{claim}{Claim}

\newtheorem*{theoremintro}{Theorem}

\theoremstyle{definition}
\newtheorem{definition}[theorem]{Definition}
\newtheorem{example}[theorem]{Example}

\newtheorem{question}[theorem]{Question}
\newtheorem*{notation}{Notation}

\theoremstyle{remark}
\newtheorem{remark}[theorem]{Remark}

\newcommand{\ydg}{{}^{\k\Gamma}_{\k\Gamma}\mathcal{YD}}

\newcommand\id{\operatorname{id}}

\newcommand\Hom{\operatorname{Hom}}
\newcommand\End{\operatorname{End}}
\newcommand\Rep{\operatorname{Rep}}

\renewcommand\c{\mathfrak{c}}

\newcommand\ufo{\mathfrak{ufo}}

\newcommand\qdim{\operatorname{qdim}}
\newcommand\qtr{\operatorname{qtr}}

\newcommand\tr{\operatorname{tr}}

\def\k{\Bbbk}

\def\ot{\otimes}

\def\s{\mathbb{S}}

\def\B{\mathfrak{B}}

\def\H{\mathbb{H}}
\def\eps{\epsilon}

\def\qb{\mathfrak{q}}
\renewcommand{\O}{\mathcal{O}}
\def\I{\mathbb{I}}

\newcommand\soc{\operatorname{soc}}

\newcommand\ibjp{\boldsymbol{|\,\rangle}}
\newcommand\bijp{\boldsymbol{\langle \,|}}
\newcommand\bibjp{\boldsymbol{\langle \,|\,\rangle}}

\newcommand\ijp{\boldsymbol{|}}

\newcommand\ibj{\text{\resizebox{\width}{.8\height}{$\ibjp$}}}
\newcommand\bij{\text{\resizebox{\width}{.8\height}{$\bijp$}}}
\newcommand\bibj{\text{\resizebox{\width}{.8\height}{$\bibjp$}}}
\renewcommand\ij{\text{\resizebox{\width}{.8\height}{$\ijp$}}}

\makeatletter
\newcommand*{\mathreflect}[1]{%
	\binrel@{#1}\binrel@@{\mathpalette\math@reflect{#1}}%
}
\newcommand*{\math@reflect}[2]{\reflectbox{\m@th$#1#2$}}
\newcommand*{\mathrotate}[3][]{%
	\binrel@{#3}\binrel@@{\vphantom{#3}\mathpalette\math@rotate{{#1}{#2}{#3}}}%
}
\newcommand*{\math@rotate}[2]{\math@@rotate#1#2}
\newcommand*{\math@@rotate}[4]{
	\sbox\z@{$\m@th#1#4$}%
	\smash{\makebox[\wd\z@]{\rotatebox[#2]{#3}{$\m@th#1#4$}}}%
}
\makeatother

\DeclareRobustCommand\longtwoheadrightarrow
{\relbar\joinrel\twoheadrightarrow}

\DeclareRobustCommand\longrightarrow
{\relbar\joinrel\rightarrow}

\newcommand{\longhookrightarrow}{\lhook\joinrel\longrightarrow}

\newcommand{\ku}{{\Bbbk}}

\newcommand{\Z}{{\mathbb Z}}

\newcommand{\N}{{\mathbb N}}

\newcommand{\Ext}{\mbox{\rm Ext\,}}

\renewcommand{\lg}{\langle}
\newcommand{\rg}{\rangle}

\def\pf{\begin{proof}}
	\def\epf{\end{proof}}

\newcommand\diag{\operatorname{diag}}

\def\ima{\mathrm{i}}

\def\G{\mathbb{G}}

\newcommand{\uRep}{\operatorname{\underline{Rep}}}
\def\g{\mathfrak{g}}

\def\ben{\begin{enumerate}[leftmargin=*]}
\def\een{\end{enumerate}}

\def\bit{\begin{itemize}[leftmargin=*]}

\def\eit{\end{itemize}}

\allowdisplaybreaks

\begin{document}
	
\title[Representations of type $A_2$]{Irreducible representations of pointed Hopf algebras of type $A_2$}

\author[Garc\'ia Iglesias, Rodriguez]{Agust\'in Garc\'ia Iglesias\,$^1$,
	Alfio Antonio Rodriguez\,$^2$} 

\address{{\em $^{1,2}$\,FaMAF-CIEM (CONICET), Universidad Nacional de C\'ordoba,
		Medina A\-llen\-de s/n, Ciudad Universitaria, 5000 C\' ordoba, Rep\'ublica Argentina.  }}

\begin{abstract}
	We classify the irreducible representations of a family of finite-dimensional pointed liftings $H_\lambda$ of the Nichols algebra associated with the diagram $A_2$ with parameter $q=-1$. 
	We show that these algebras have infinite representation type and construct an indecomposable $H_\lambda$-module of dimension $n$ for each $n\in\N$.
	Finally, we study a semisimple category $\uRep H_\lambda$ arising as a quotient of $\Rep H_\lambda$.
\end{abstract}


	\maketitle

	\section{Introduction}\label{sec:intro}

	For $N,M\in\N$, we study the representations of pointed Hopf algebras over 
	\[
	\Gamma=\lg g_1,g_2:g_1g_2=g_2g_1, \ g_1^{2N}=g_2^{2M}=1\rg\simeq \Z/2N\Z\times \Z/2M\Z
	\]
	and with infinitesimal braiding 
	of diagonal type $A_2$ with parameter $q=-1$.
	
	These are deformations of the positive part  of the small quantum group $\mathfrak{u}_{\sqrt{-1}}(\mathfrak{sl}_3)$ and are classified in terms of triples $\lambda=(\lambda_1,\lambda_2,\lambda_3)\in \k^3$: they are the quotients of the smash product $\k\lg a_1,a_2\rg\#\k\Gamma$ with relations:
	\begin{align}\label{eqn:rels_conm}
		g_1a_1&=-a_1g_1, & g_1a_2&=-a_2g_1, &  g_2a_1&=a_1g_2, &  g_2a_2&=-a_2g_2,
	\end{align}
	by the  relations, see \S\ref{sec:Hlambda} for details:
	\begin{align}\label{eqn:rels_def}
		\begin{split}
			a_1^{2}=\lambda_1(1-&g_1^2), \qquad  a_2^2=\lambda_2(1-g_2^2),\\
			a_1a_2a_1a_2+a_2a_1a_2a_1&=\lambda_3(1-g_1^2g_2^2)-2\lambda_1\lambda_2(1+g_2^2)(1-g_1^2).
		\end{split}
	\end{align}
	Let $\G_{2N}, \G_{2M}\subset\k$ stand for the groups of $2N$-th and $2M$-th roots of unity; we refer to \S\ref{sec:graphics} for the graphical notation~\eqref{eqn:simpledraw} for $H_\lambda$-modules. A simplified version of our classification result, Theorem \ref{thm:simples}, reads as follows:
	\begin{theoremintro}
		Define $\O=\G_{2N}\times \G_{2M}$. For each  $\lambda\in\k^3$ there is a decomposition $\O=\O_1\sqcup \O_2\sqcup \O_4$ so that, up to isomorphism, the simple modules of $H_\lambda$ 
		are $|\O_1|$ modules $L_1^\chi$ of dimension 1, $\frac{1}{2}|\O_2|$ modules $L_{2,h}^\chi$ or $L_{2,v}^\chi$ of dimension 2 and $\frac{1}{2}|\O_4|$ modules $L_4^\chi(d)$ of dimension 4 (or $\frac{1}{2}|\O_4|- 2$ in the exceptional case of~\eqref{eqn:number=4}).
		These modules are represented as 	
		\begin{align}\label{eqn:simpledraw}
			L_1^\chi: & \, \ij
			&
			L_{2, h}^{\chi}: & \xymatrix{
				\ij\ar@{<->}[r]_{\alpha_1}&
				\bij}  \quad 	 	
			L_{2, v}^{\chi}: \,  \xymatrix{
				\ij\ar@{<->}[d]_{\alpha_2}\\
				\bibj
			} 
			&
			L_{4}^{\chi}(d):& 
			\xymatrix{
				\ij\ar@{<->}[r]_{\alpha_1} \ar@{<->}[d]_{c}^{d}&\bij\ar@{<->}[d]_{\alpha_2}\ar@{}[dl]\\
				\bibj\ar@{<->}[r]^{\alpha_1} &\ibj
			}
		\end{align}
		where $\alpha_1,\alpha_2,c,d \in \k$ depend on  $(\lambda,\chi)$; see~\eqref{eqn:alpha} and~\eqref{eqn:dim4-quadratic}.
	\end{theoremintro}
	\subsection{Organization}
	In \S\ref{sec:simples}, we define the simple modules of dimensions 1, 2, and 4. The classification is proved in \S\ref{sec:classif}, using projective covers. In \S\ref{sec:gabriel}, we compute the Gabriel quiver and show that these algebras are not of finite representation type. We study indecomposable modules in \S\ref{sec:indecomposables} and define an indecomposable $H_\lambda$-module of dimension $n$ for each $n\geq 1$. We classify indecomposable modules of small dimension. We show that $H_\lambda$ is spherical when $N$ is odd; it thus gives rise to a semisimple category $\uRep H_\lambda$, see \S\ref{sec:quotient}. We compute part of its fusion rules using results from the previous section.
	
	\subsection{Background}
	
	In~\cite{AB} the authors study representations of liftings of quantum planes; here the Dynkin diagram is a finite product $A_1\times\dots\times A_1$. 
	The representation theory of a large class of pointed liftings of diagonal type is the subject of~\cite{ARS}.
	
	Some work along these lines has also been carried out for related Nichols algebras and their Drinfeld doubles.
	The simple modules of the Drinfeld double of the Jordan plane are classified in~\cite{ADP}, while in~\cite{AP} the authors classify an infinite family of indecomposable modules over this algebra.
	Analogous results for the super Jordan plane are obtained in~\cite{ABFD}, and the simple modules of the Drinfeld double of the Nichols algebra of type $\ufo(7)$ are computed in~\cite{AAMR}.

	In~\cite{GI,GR} we followed these ideas to analyze the representation theory for liftings of the Fomin--Kirillov algebra on three generators. In this case the braiding is non-diagonal and the underlying group is not abelian; it projects onto the symmetric group $\s_3$. The liftings are Hopf cocycle deformations of the graded algebra associated to the coradical filtration: in loc. cit. we found a connection between the number of simple modules and the expression of the cocycle as an exponential of a Hochschild 2-cocycle. We continue this analysis in Corollary \ref{cor:cocycles}.
	
	When $\g$ is a finite-dimensional simple Lie algebra, there is a vast collection of results for $\Rep U_q(\g)$ and $\Rep u_q(\g)$ and their connection to the representation theory of the corresponding algebraic group, or the representations of $\g$ in positive characteristic, see e.g.~\cite{dl,l1,rosso}.

	\section{Preliminaries}
	
	We work over an algebraically closed field $\k$	of characteristic zero.
	We recall a general result, including a short proof for completeness.
	
	\begin{lemma}\label{lem:quotients}
		Let $A$ be a $\k$-algebra. Assume there is a finite group $G$ such that $\k G\subset A$ as a subalgebra. 
		If $L$ is an irreducible $A$-module, then there is a simple module $S\in\widehat{G}$ such that $L$ is a quotient of the induced $A$-module $A(S)\coloneqq {}_AA\ot_{\ku G} S$.
	\end{lemma}
	\pf
	Let $L$ be an irreducible $A$-module, and consider its decomposition as a $G$-module: $L_{|}\simeq S_1\oplus \dots \oplus S_k$, with each $S_i\in\widehat{G}$.
	The projection $A\ot_{\ku G}L=\bigoplus_{i=1}^k A(S_i)\twoheadrightarrow L$ induces morphisms $A(S_i)\hookrightarrow  A\ot_{\ku G}L \twoheadrightarrow L$. 
	These maps cannot all be zero, so there exists some index $j$ such that $A(S_j)\twoheadrightarrow L$. 
	\epf
	
	\subsection{The algebras $H_\lambda$}\label{sec:Hlambda} 
	Finite-dimensional, non-semisimple pointed Hopf algebras with a fixed abelian group of group-like elements $\Gamma$ are classified in terms of (infinitesimal) braiding matrices $\qb$ (or, equivalently,  generalized Dynkin diagrams). They are liftings of the corresponding Nichols algebra $\B_{\qb}$, with a realization $\B_{\qb}\in\ydg$.
	We refer the reader to~\cite{AS,An,AnG,H} for further details.

	In this article we study the irreducible representations of pointed Hopf algebras $H$ whose group of group-like elements is $\Gamma\coloneqq \Z/2N\Z\times\Z/2M\Z$, $N,M\geq 1$, and whose infinitesimal braiding is the matrix $\qb=\begin{psmallmatrix}
		q_{11}&q_{12}\\q_{21}&q_{22}
	\end{psmallmatrix}\in\k^{2\times 2}$  
	associated to the diagram
	\[
	\Dchaintwo{-1}{-1}{-1}
	\]
	That is, $q_{11}=q_{22}=-1=q_{12}q_{21}$. The corresponding Nichols algebra is 
	\[
	\B_{\qb}=\lg x_1,x_2|x_1^2=x_2^2=x_1x_2x_1x_2+q_{12}^2x_2x_1x_2x_1=0\rg.
	\]
	This algebra has dimension 8. A linear basis is given by the set 
	\begin{align}\label{eqn:basis}
		\mathbb{B}&=
		\{1,x_1,x_2,x_1x_2,x_2x_1,x_1x_2x_1,x_2x_1x_2,x_2x_1x_2x_1\}.
	\end{align}
	We remark that $\B_{\qb}$ is the positive part $\mathfrak{u}^+_{\sqrt{-1}}(\mathfrak{sl}_3)$ of the small quantum group $\mathfrak{u}_{\sqrt{-1}}(\mathfrak{sl}_3)$.
	The basis $\mathbb{B}$ coincides with the usual PBW basis in this context.
	
	When $q_{12}\neq\pm1$, $\B_{\qb}$ admits no deformations, namely $\B_{\qb}\#\k\Gamma$ is, up to isomorphism, unique in this class. The same holds when $N=M=1$~\cite{AD}.

	\subsubsection*{Convention}\label{page:convention}
	We shall fix $q_{12}=-1, q_{21}=1$ and assume $N,M\geq 2$.
	
	\medbreak
	
	The liftings of $\B_{\qb}$ over $\Gamma$ are classified by triples $\lambda=(\lambda_1,\lambda_2,\lambda_3)\in\k^3$: these are the algebras $H_\lambda$ as in~\eqref{eqn:rels_conm} and~\eqref{eqn:rels_def}. The  coalgebra structure is determined by, for $i=1,2$:
	\begin{align}\label{eqn:comultiplication}
		\Delta(g_i)&=g_i\ot g_i, & \Delta(a_i)=a_i\ot 1+g_i\ot a_i.
	\end{align}
	Therefore, $\eps(g_i)=1, S(g_i)=g_i^{-1}$,  and $\eps(a_i)=0$, $S(a_i)=-g_i^{-1}a_i$.
	
	\begin{remark}
		The symmetric case $q_{12}=1, q_{21}=-1$ is equivalent under the exchange $1\leftrightsquigarrow 2$ and $N\leftrightsquigarrow M$. 
		When $N=1$, we may set $\lambda_1=0$; similarly, when $M=1$, we may set $\lambda_2=0$.
	\end{remark}

	\subsubsection{Cocycle deformations}\label{sec:cocycles}
	
	We recall that the algebras $H_\lambda$ are Hopf cocycle deformations of the graded Hopf algebra $B=\B_{\qb}\#\k\Gamma$; see~\cite{GS} for details and background. We point out in Corollary \ref{cor:cocycles} a connection between the \emph{shape} of the cocycle and the representation theory of the algebra.
	
	We briefly recall this setting. Let $H$ be a Hopf algebra: if $\eta$ is a Hochschild 2-cocycle, then under certain conditions the exponential $e^\eta\colon H\ot H\to \k$ is a Hopf cocycle. This technique is far from exhausting the set of cocycles for $H$: a Hopf 2-cocycle not cohomologous to an exponential is called pure. 
	
	The Hopf cocycles for $H_\lambda$ were classified in~\cite{GS}: it is shown in \emph{loc. cit.} that $H_\lambda$ can be obtained from $B$ with an exponential Hopf 2-cocycle if and only if at most one of the parameters $\lambda_i$ (for $i = 1, 2, 3$) is nonzero, see~\cite[Theorem 6.13]{GS}. Otherwise, the cocycle deforming $B$ into $H_\lambda$ is pure.

	\subsection{Simple $\Gamma$-modules}
	We set $\I_\theta=\{1,\dots,\theta\}\subset \N$, $\I_\theta^\circ=\I_\theta\cup\{0\}$, $\theta\in\N$.
	
	Let $\zeta, \xi$ be primitive roots of 1 of orders $2N$ and $2M$, respectively, so $\zeta^{2N}=\xi^{2M}=1$.
	As $\Gamma$ is abelian, every simple module is one-dimensional and these are parametrized by the set 
	\begin{align*}
		\O=\{(\zeta^i,\xi^j):(i,j)\in \I_{2N-1}^\circ\times \I_{2M-1}^\circ\}=\G_{2N}\times \G_{2M}.
	\end{align*} 
	If $\chi=(\zeta^i,\xi^j)\in\O$, then the simple module $S=S_\chi$ can be described via:
	\begin{align}\label{eqn:Schi}
		S_\chi=\k\lg z_{ij}\rg, \qquad g_1\cdot z_{ij}=\zeta^i z_{ij}, \quad g_2\cdot z_{ij}=\xi^j z_{ij}.
	\end{align}

	\section{Simple modules of small dimension}\label{sec:simples}
	
	In this section we study irreducible modules of small dimension. These will account for all irreducible $H_\lambda$-modules.
	We introduce a collection of scalars. 
	For each $\chi=(\zeta^i,\xi^j)\in\O$, we define $\alpha_t=\alpha_t(\chi)\in\k$, $t\in\I_3$, as
	\begin{align}\begin{split}\label{eqn:alpha}
			&\alpha_1=\lambda_1(1-\zeta^{2i}),\quad \alpha_2=\lambda_2(1-\xi^{2j}),\\ &\alpha_3=\lambda_3(1-\zeta^{2i}\xi^{2j})-2\lambda_1\lambda_2(1+\xi^{2j})(1-\zeta^{2i}).	
		\end{split}	
	\end{align}
	
	\subsection{Isotypic components}
	For any $H_\lambda$-module $M$, we have the decomposition into $\Gamma$-isotypic components $M_{|}=\bigoplus_{{\chi\in\O}}M_{|}[\chi]$; here 
	if $\chi=(\chi_1,\chi_2)\in\O$, then $g_i$ acts on $M_{|}[\chi]$ as multiplication by $\chi_i$, $i\in\I_2$. 
	We set  $\overline{\chi}\coloneqq (-\chi_1, \chi_2)$.
	
	The generators $a_1, a_2$ act by shifting components as follows: $a_1\cdot M_{|}[\chi]\subset M_{|}[\overline\chi]$ and $a_{2}\cdot M_{|}[\chi]\subset M_{|}[-\chi]$. That is, we have the following interaction between $\Gamma$-isotypic components:
	\begin{align}\begin{split}\label{eqn:squares}
			\xymatrix{
				\chi\ar@{<.>}[r]^{a_1} \ar@{<.>}[d]_{a_2}&\overline{\chi}\ar@{<.>}[d]^{a_2}\\
				-\chi\ar@{<.>}[r]_{a_1} & -\overline{\chi}
			}
	\end{split}\end{align}
	
	\begin{notation}
		If $\chi\in\O$, we let $\Omega(\chi)=\{\chi,\overline{\chi},-\overline{\chi},-\chi\}$.
		We denote by $M^\chi$ the (possibly zero) submodule generated by the components $M_{|}[\chi']$, $\chi'\in\Omega(\chi)$, so
		$M^\chi_{|}=M_{|}[\chi]\oplus M_{|}[\overline{\chi}]\oplus M_{|}[-\overline{\chi}]\oplus M_{|}[-\chi]$.	Note that $M^\chi$ can still be decomposable.
		We write $\widehat{\O}$ for a set of representatives of the classes of the equivalence relation on $\O$ given by $\chi\sim \chi'$ if and only if $\chi'\in\Omega(\chi)$.
	\end{notation}
	
	\begin{lemma}\label{lem:indecomp}
		Let $M$ be an $H_\lambda$-module. Then $M\simeq \bigoplus_{\chi\in \widehat{\O}} M^{\chi}$.
		Hence, if $M$ is indecomposable, then there exists $\chi\in\O$ such that $M=M^\chi$. \hfill\qed
	\end{lemma}
	
	\begin{remark}\label{rem:dim-Hchi}
		This decomposition implies that
		$
		H_\lambda^\chi=\bigoplus\limits_{\mathclap{\chi'\in\Omega(\chi)}} H_\lambda\ot_{\Gamma} S_{\chi'}.
		$
		In particular, $\dim H_\lambda^\chi=32$ for any $\chi\in\O$. 
		Furthermore, this lemma allows us to specify the decomposition 
		$
		H_\lambda=\bigoplus_{L\,\text{irr.}}P(L)^{\oplus\dim L},
		$ 
		where $P(L)$ stands for the projective cover of the irreducible $H_\lambda$-module $L$. 
		
		Our goal is to refine the identity
		$\dim H_\lambda=\sum_{L\,\text{irr.}}\dim L\dim P(L)$. 
		By the lemma, for each $L$ and its projective cover $P(L)$, there is $\chi$ such that $L=L^\chi$ and $P(L)=P(L)^\chi$. 
		We can restrict this decomposition to obtain
		\begin{equation}\label{eqn:dim32}
			H_\lambda^\chi=\bigoplus
			P(L)^{\oplus\dim L}, \qquad \text{and hence} \qquad 32=\sum
			\limits
			\dim L\dim P(L),
		\end{equation}
		where the sums run over all irreducible modules 
		$L$ with support in $\Omega(\chi)$.
	\end{remark}

	\subsection{One-dimensional modules}
	We begin by classifying the simple modules of dimension 1. This classification is straightforward.

	\begin{notation} We set:
		\begin{align}\label{eqn:O1}
			\O_1\coloneqq\{\chi\in\O|\alpha_1=\alpha_2=\alpha_3=0\}.
		\end{align}
	\end{notation}

	\begin{proposition}\label{pro:dim1}
		Let $L$ be a 1-dimensional $H_\lambda$-module. Then there exists $\chi \in \O$ such that $L_{|} \simeq S_\chi$ and both $a_1$ and $a_2$ act trivially on $L$. We denote this module by $L^\chi_{1}$.
		Such a module $L_1^\chi$ exists if and only if $\alpha_1=\alpha_2=\alpha_3=0$, that is, if and only if $\chi\in\O_1$.
		We have $L_1^\chi\simeq L_1^{\chi'}$ if and only if $\chi=\chi'$.
	\end{proposition}
	\pf
	The first part follows from the description of the $\Gamma$-modules and the observation in~\eqref{eqn:squares}.
	If $L=\lg z\rg$, then  $\alpha_1(\chi)=\alpha_2(\chi)=\alpha_3(\chi)=0$, by~\eqref{eqn:rels_def}. 
	\epf

	\subsection{Two-dimensional modules}
	We consider modules of dimension 2.
	
	\begin{proposition}\label{pro:dim2}
		Let $L$ be a simple $H_\lambda$-module of dimension 2. Then there exists $\chi \in \O$ such that $\alpha_3 = 0$ and one of the following holds:
		\begin{enumerate}[leftmargin=*]
			\item[(a)] $\alpha_1\neq 0$, $\alpha_2=0$ and $L$  has a basis $\{v,w\}$ with $\lg v\rg_{|}\simeq S_\chi$, $\lg w\rg_{|}\simeq S_{\overline{\chi}}$ and
			\begin{align}\label{eqn:L2h}
				a_1\cdot v&=w, &  a_2\cdot v&=0,&   a_1\cdot w&=\alpha_1v, & a_2\cdot w&=0.
			\end{align}
			We denote this $L$ by $L_{2,h}^\chi$. Moreover, $L_{2,h}^\chi\simeq L_{2,h}^{\phi}$ if and only if $\phi\in\{\chi,\overline{\chi}\}$.
			\item[(b)] $\alpha_2\neq 0$, $\alpha_1=0$ and $L$  has a basis $\{v,w\}$ with $\lg v\rg_{|}\simeq S_\chi$,  $\lg w\rg_{|}\simeq S_{-\chi}$  and
			\begin{align}\label{eqn:L2v}
				a_1\cdot v&=0, &  a_2\cdot v&=w, & a_1\cdot w&=0, & a_2\cdot w&=\alpha_2 v.
			\end{align}
			We denote this $L$ by $L_{2,v}^\chi$. Moreover, $L_{2,v}^\chi\simeq L_{2,v}^{\phi}$ if and only if $\phi\in\{\chi,-\chi\}$. 
		\end{enumerate} 
	\end{proposition}
	\pf
	It is easy to check that, under the preceding hypotheses for each case, the assignments in~\eqref{eqn:L2h} and~\eqref{eqn:L2v} yield 2-dimensional simple modules.
	
	For the converse, let $0\neq v\in L$. We may assume that there is $\chi\in\O$ such that $\lg v\rg_{|}\simeq S_\chi$. As $L$ is simple of dimension two, we get that either $a_1\cdot v\neq 0$ or $a_2\cdot v\neq 0$. Moreover, only one of them is nonzero, since by~\eqref{eqn:squares} the nonzero vectors among $v$, $a_1\cdot v$, $a_2\cdot v$ lie in different isotypic components, and hence are linearly independent.
	
	Say $a_1\cdot v\neq 0$ and set $w\coloneqq a_1\cdot v$; thus $\{v,w\}$ is a basis of $L$, where $\lg w\rg_{|}\simeq S_{\overline{\chi}}$ by~\eqref{eqn:squares}. In particular, $\alpha_1\neq 0$: indeed, $a_1\cdot w\neq 0$, since otherwise $\lg w\rg\subset L$ would be a submodule, and $a_1\cdot w=a_1^2\cdot v=\alpha_1v$. The case $a_2\cdot v\neq 0$ is analogous, and gives rise to $L_{2,v}^\chi$.  
	
	In either case, as $0=(a_1a_2)^2+(a_2a_1)^2$ in this module, we have $\alpha_3=0$.
	The remaining assertions follow directly from the definitions.
	\epf
	
	\begin{notation} We set:
		\begin{align}\label{eqn:O2}
			\O_2\coloneqq\{\chi\in\O| \alpha_3=0 \text{ and } (\alpha_1\neq 0=\alpha_2 \text{ or } \alpha_2\neq 0=\alpha_1)\}.
		\end{align}
		We write $\overline{\O}_2$ for a subset of representatives of isomorphism classes of 2-dimensional simple modules, so $|\overline{\O}_2|=\frac{1}{2}|\O_2|$ by Proposition \ref{pro:dim2}.
	\end{notation}

	\subsection{Four-dimensional modules}
	
	We now focus on those $\chi \in \O$ such that $\chi \notin \O_1 \sqcup \O_2$; this defines the subset $\O_4 \subset \O$, given by
	\begin{align}\label{eqn:O4}
		\O_4=\{\chi\in\O: \text{either }&\alpha_1\alpha_2\alpha_3\neq 0\text{ or }\alpha_1=\alpha_2=0\text{ and }\alpha_3\neq 0 \\ 
		\notag &\text{ or at most one }\alpha_t, t\in\I_3, \text{ is zero}\}.
	\end{align}

	\begin{proposition}\label{pro:dim4}
		Let $L$ be an irreducible $H_\lambda$-module of dimension 4. 
		Then there exist $\chi \in \O_4$ and $(d, c) \in \k^2$ such that
		\begin{align}\label{eqn:dim4-quadratic}
			\alpha_1^2\alpha_2d^2-\alpha_3d+\alpha_2&=0, & c&=\alpha_3-\alpha_1^2\alpha_2d,
		\end{align}
		so that $L$ has a basis $\{v_1,v_2,v_3,v_4\}$ with 
		\begin{align*}
			\lg v_1\rg_{|}&\simeq S_\chi, & \lg v_2\rg_{|}&\simeq S_{\overline{\chi}}, & \lg v_3\rg_{|}&\simeq S_{-\overline{\chi}}, & \lg v_4\rg_{|}&\simeq S_{-\chi}, 
		\end{align*}
		and such that 
		\begin{align}\label{eqn:actiond}
			a_1\cdot v_1=v_2, \ a_2\cdot v_2=v_3, \ a_1\cdot v_3=v_4, \ a_2\cdot v_1=d\,v_4, \ a_2\cdot v_4=c\, v_1.
		\end{align}
		Conversely, given $\chi\in\O_4$ and $d$ as in~\eqref{eqn:dim4-quadratic}, the equations above define an irreducible $H_\lambda$-module $L_4^\chi(d)$ with basis  $\{v_1,v_2,v_3,v_4\}$. 
	\end{proposition}
	We shall look into isomorphism classes in Proposition \ref{pro:iso4}.
	\pf
	Let $L$ be such a module. Then there is $\psi\in\O$ for which $L_{|}[\psi]\neq 0$. Fix $0\neq v_1\in L_{|}[\psi]$. Observe that we cannot have $a_1\cdot v_1=0$ and $a_2\cdot v_1=0$, as otherwise $L_1^\psi\simeq \lg v_1\rg$ is a submodule. Also, notice that $L_{|}\simeq S_\psi\oplus S_{\overline{\psi}}\oplus S_{-\overline{\psi}}\oplus S_{-\psi}$, by~\eqref{eqn:squares}.
	
	Assume $a_1\cdot v_1\neq 0$. Then $v_2\coloneqq a_1\cdot v_1$ is such that $\lg v_2\rg_{|}\simeq S_{\overline{\psi}}$. Now we may assume that $a_2\cdot v_2\neq 0$. Indeed, if $a_2\cdot v_2=0$, then $\lg v_2\rg$ is not a submodule, so $0\neq a_1\cdot v_2=a_1^2\cdot v_1=\alpha_1 v_1$, and $\lg v_1,v_2\rg$ is not a submodule, so $a_2\cdot v_1\neq 0$. Hence we can replace $\psi$ by $\overline{\psi}$ and $v_1$ by $v_2$, as $a_1\cdot v_2\neq 0$ and $a_2a_1\cdot v_2=\alpha_1\,a_2\cdot v_1\neq 0$. Set $v_3\coloneqq a_2\cdot v_2$, so $\lg v_3\rg_{|}\simeq S_{-\overline{\psi}}$. A similar argument shows that $v_4\coloneqq a_1\cdot v_3\neq 0$, $\lg v_4\rg_{|}\simeq S_{-\psi}$ and $a_2\cdot v_4\in\lg v_1\rg$ by a dimension argument. Similarly, $a_2\cdot v_1\in\lg v_4\rg$.
	
	Fix $c, d \in \k$ such that $a_2\cdot v_4=cv_1$ and $a_2\cdot v_1=dv_4$. Thus, we necessarily have $cd=\alpha_2$. Hence 
	the actions of $a_1$ and $a_2$ are given by the matrices
	\begin{align}\label{eqn:a1a2-dim4}
		[a_1]&=A\coloneqq\begin{psmallmatrix}
			0&\alpha_1&0&0\\
			1&0&0&0\\
			0&0&0&\alpha_1\\
			0&0&1&0
		\end{psmallmatrix}, & 
		[a_2]&=B\coloneqq\begin{psmallmatrix}
			0&0&0&c\\
			0&0&\alpha_2&0\\
			0&1&0&0\\
			d&0&0&0
		\end{psmallmatrix}.
	\end{align}
	Now the relation $(a_1a_2)^2+(a_2a_1)^2=\alpha_3$ gives $c+\alpha_1^2\alpha_2d=\alpha_3$. Thus for each solution $d\in\k$ of~\eqref{eqn:dim4-quadratic}
	the matrices~\eqref{eqn:a1a2-dim4} with $c=\alpha_3-\alpha_1^2\alpha_2d$ determine the module $L$. Note that it is necessary that either $\alpha_1 = \alpha_2 = 0$ and $\alpha_3 \ne 0$, or at most one of the parameters $\alpha_t$, $t\in\I_3$, vanishes.
	Indeed:
	\begin{itemize}[leftmargin=*]
		\item If $\alpha_1=\alpha_2=\alpha_3=0$, then $L$ is not simple as $\lg v_4\rg$ is a submodule.
		\item If $\alpha_1=\alpha_3=0$, then~\eqref{eqn:dim4-quadratic} implies $\alpha_2=0$; hence $L$ is not simple.
		\item If $\alpha_2=\alpha_3=0$, then $L$ is not simple, as $\lg v_3, v_4\rg$ is a submodule.
	\end{itemize}
	Otherwise, the module is simple. Thus the statement of the proposition follows in this case, for $\chi\coloneqq \psi\in\O_4$.
	
	Now, assume that $a_1\cdot v_1=0$; hence $\alpha_1=0$, and thus $\alpha_3\neq 0$, $d=\alpha_2/\alpha_3$.
	
	Setting $w_1=v_1$, $w_4=a_2\cdot w_1$, $w_3=a_1\cdot w_4$ and $w_2=a_2\cdot w_3$, we obtain as above a linearly independent set $\{w_1,w_2,w_3,w_4\}$ with $\lg w_1\rg_{|}\simeq S_\psi$, $\lg w_2\rg_{|}\simeq S_{\overline{\psi}}$, $\lg w_3\rg_{|}\simeq S_{-\overline{\psi}}$, $\lg w_4\rg_{|}\simeq S_{-\psi}$ and such that the action is encoded by the matrices
	\begin{align*}
		[a_1]&=\begin{psmallmatrix}
			0&\alpha_3&0&0\\
			0&0&0&0\\
			0&0&0&1\\
			0&0&0&0
		\end{psmallmatrix}, & 
		[a_2]&=\begin{psmallmatrix}
			0&0&0&\alpha_2\\
			0&0&1&0\\
			0&\alpha_2&0&0\\
			1&0&0&0
		\end{psmallmatrix}.
	\end{align*}
	Let us set $v_1\coloneqq w_4$, $v_2\coloneqq w_3$, $v_3\coloneqq w_2$, $v_4\coloneqq \alpha_3 w_1$.
	This is a new basis for which $a_1\cdot v_1=v_2$, $a_2\cdot v_2=v_3$, $a_1\cdot v_3=v_4$ and
	\begin{align*}
		a_2\cdot v_1&=a_2\cdot w_4=a_2^2\cdot w_1=\alpha_2 w_1=(\alpha_2/\alpha_3) v_4=d v_4,\\
		a_2\cdot v_4&=\alpha_3 a_2\cdot w_1=\alpha_3 w_4=\alpha_3 v_1=(\alpha_3-\alpha_1^2\alpha_2d)v_1.
	\end{align*}
	Hence~\eqref{eqn:actiond} also holds, by setting $\chi\coloneqq -\psi$. 
	
	The converse is straightforward, by checking the relations. 
	\epf
	
	Next we study the isomorphism classes of the modules $L_4^\chi(d)$, $\chi\in\O_4$. 
	
	We fix some notation: if $\chi=(\zeta^i,\xi^j)\in\O_4$ and $\zeta^{2i}\xi^{2j}\neq 1$, then we set
	\begin{align}\label{eqn:chipm}
		\chi_+&\coloneqq\tfrac{1-\zeta^{2i}}{1-\zeta^{2i}\xi^{2j}}, & \chi_-&\coloneqq\xi^{2j}\tfrac{1-\zeta^{2i}}{1-\zeta^{2i}\xi^{2j}}.
	\end{align}
	
	\begin{proposition}\label{pro:iso4}
		Let $\chi,\phi\in \O_4$, and let $d=d(\chi), e=e(\phi)$ be as in~\eqref{eqn:dim4-quadratic}. If $L_4^\chi(d)\simeq L_4^\phi(e)$, then $\phi\in\Omega(\chi)$.
		Moreover,
		\begin{itemize}
			\item[(i)] $L_4^\chi(d)\simeq L_4^\chi(e)$ if and only if $d=e$.
			\item[(ii)] $L_4^\chi(d)\simeq L_4^{-\overline{\chi}}(e)$ if and only if $d=e$.
			\item[(iii)] $L_4^\chi(d)\simeq L_4^{-\chi}(e)$ if and only if $\alpha_1\alpha_2\neq 0$ and $\alpha_1^2de=1$.
			\item[(iv)] $L_4^\chi(d)\simeq L_4^{\overline{\chi}}(e)$ if and only if $\alpha_1\alpha_2\neq 0$ and $\alpha_1^2de=1$.
		\end{itemize}	
		Therefore, if $\overline{\O}_4$ stands for a set of representatives of isomorphism classes of irreducible $H_\lambda$-modules of dimension 4, then	
		\begin{equation}\label{eqn:number=4}
			|\overline{\O}_4|=\begin{cases}
				\frac{1}{2}|\O_4|-2, & \text{if } \lambda_3=4\psi_\pm\lambda_1\lambda_2 \text{ for some } \psi\in\O_4,\\
				\frac{1}{2}|\O_4|, & \text{otherwise.}
			\end{cases}
		\end{equation}
	\end{proposition}
	\pf
	The first assertion follows from the decomposition into $\Gamma$-components.
	
	The  modules $L_4^\chi(d)$ are generated by a vector $v_1$ such that $\lg v_1\rg_{|}\simeq S_\chi$ as $\Gamma$-modules. Hence any morphism  $f\colon L_4^\chi(d)\to L_4^{\chi}(e)$ is completely determined by $f(v_1)$, as $f(v_2)=a_1\cdot f(v_1)$, $f(v_3)=a_2\cdot f(v_2)$ and $f(v_4)=a_1\cdot f(v_3)$. This implies (i). Case (ii) follows analogously: if $\{w_1,\dots,w_4\}$ is the basis for $L_4^{-\overline{\chi}}(e)$, then $\lg w_3\rg_{|}\simeq S_\chi$ and we can assume $f(v_1)=w_3$.

	The case in (iii) is similar. If $f\colon L_4^\chi(d)\to  L_4^{-{\chi}}(e)$ and  $\{w_1,\dots, w_4\}$ is the basis for $L_4^{-{\chi}}(e)$, then $\lg w_4\rg_{|}\simeq S_\chi$ and hence we can assume $f(v_1)=w_4$. Therefore $f(v_2)=a_1\cdot w_4=\alpha_1 w_3$, $f(v_3)=\alpha_1 a_2\cdot  w_3=\alpha_1\alpha_2w_2$ and $f(v_4)=\alpha_1^2\alpha_2w_1$. In particular, $\alpha_1\alpha_2\neq 0$. 
	
	Now, $f(a_2\cdot v_1)=df(v_4)=d\alpha_1^2\alpha_2w_1$ and $a_2\cdot f(v_1)=a_2\cdot w_4=c'w_1$, while $f(a_2\cdot v_4)=cf(v_1)=cw_4$ and $a_2\cdot f(v_4)=\alpha_1^2\alpha_2 a_2\cdot w_1=\alpha_1^2\alpha_2 ew_4$. Here $c'=\alpha_3-\alpha_1^2\alpha_2e$. This gives $d\alpha_1^2\alpha_2=c'$ and $c=\alpha_1^2\alpha_2 e$. Hence $e=1/(\alpha_1^2d)$.
	
	Case (iv) is analogous; here $f(v_1)=w_2$.
	
	Let us now look at isomorphism classes. First, if $\alpha_1\alpha_2=0$, then $\alpha_3\neq 0$ and equation~\eqref{eqn:dim4-quadratic} becomes $\alpha_3d-\alpha_2=0$; it thus has a single solution $d=\alpha_2/\alpha_3$. Following \emph{(i)}--\emph{(iv)}, we have two non-isomorphic 4-dimensional modules, namely $L_4^\chi(d)\ (\simeq L_4^{-\overline{\chi}}(d))$ and $L_4^{-\chi}(d)\ (\simeq L_4^{\overline{\chi}}(d))$.
	
	Assume next that $\alpha_1\alpha_2\neq 0$ and let $d,e$ be the two (possibly equal) solutions of~\eqref{eqn:dim4-quadratic}: observe that $e=\frac{1}{\alpha_1^2d}$ and let $D(\chi)$ be the discriminant:
	\begin{align}\label{eqn:discriminant}
		D(\chi)=\alpha_3^2-4\alpha_1^2\alpha_2^2.
	\end{align}
	
	If $D(\chi)\neq0$, then $d\neq e$ and we have again two isomorphism classes, given by $L_4^\chi(d)\simeq L_4^{-\overline{\chi}}(d)$ and $L_4^\chi(e)\simeq L_4^{-\chi}(d)\ (\simeq L_4^{\overline{\chi}}(d))$. 
	
	On the other hand, if $D(\chi)=0$, then $L_4^\chi(d)\simeq L_4^{-\overline{\chi}}(d)$ by \emph{(ii)} and $L_4^{-\chi}(d)\stackrel{\emph{(iii)}}{\simeq} L_4^{\chi}(d)\stackrel{\emph{(iv)}}{\simeq} L_4^{\overline{\chi}}(d)$, which determines a single 4-dimensional module up to isomorphism associated to $\chi$; here $d=2\alpha_2/\alpha_3$. Let $Z\coloneqq\{\rho\in\O_4:D(\rho)=0\}$. We state the following claim and prove it below.
	\begin{claim}
		If $\lambda_3=4\psi_\pm\lambda_1\lambda_2$ for some $\psi\in\O_4$, then $Z=\Omega(\psi)\cup\Omega(\psi^{-1})$, where $\psi^{-1}=(\psi_1^{-1},\psi_2^{-1})$; otherwise $Z=\emptyset$.
	\end{claim}
	Assume the claim holds, and let $\approx$ be the equivalence relation on $\O_4$ given by $\rho\approx\rho'$ if and only if $\rho'\in\{\rho,-\overline{\rho}\}$.
	If $\lambda_3\neq4\psi_\pm\lambda_1\lambda_2$ for all $\psi\in\O_4$, then $\overline{\O}_4$ is in bijection with $\O_4/{\approx}$.
	If, instead, $\lambda_3=4\psi_\pm\lambda_1\lambda_2$ for some $\psi\in\O_4$, then $\overline{\O}_4$ is in bijection with $(\O_4\setminus Z)/{\approx}$ together with the two classes $[L_4^{\psi}(d)]$ and $[L_4^{\psi^{-1}}(d')]$, where $d$ and $d'$ are the double roots of~\eqref{eqn:dim4-quadratic} for $\psi$ and $\psi^{-1}$, respectively. This gives formula~\eqref{eqn:number=4}.
	\medbreak
	\noindent
	\emph{Proof of the Claim.} If $\psi=(\zeta^i,\xi^j)\in Z$, then necessarily $\alpha_1\alpha_2\alpha_3\neq 0$. In particular, $\lambda_1\lambda_2\neq 0$ and $\zeta^{2i},\xi^{2j}\neq 1$.
	Here $D(\psi)=0$ becomes $\alpha_3=\pm2\alpha_1\alpha_2$, or, more precisely:
	\begin{align}\label{eqn:lambdas}
		\lambda_3(1-\zeta^{2i}\xi^{2j})=2\lambda_1\lambda_2(1-\zeta^{2i})[(1+\xi^{2j})\pm (1-\xi^{2j})].
	\end{align}
	In particular, $\lambda_3\neq 0$ and $\zeta^{2i}\xi^{2j}\neq 1$, so $\psi\in Z$ if and only if $\lambda_3=4\psi_\pm\lambda_1\lambda_2$. Thus $Z\neq \emptyset$ if and only if $\lambda_3=4\psi_\pm\lambda_1\lambda_2$ for some $\psi\in\O_4$.
	
	Assume that $\lambda_3=4\psi_+\lambda_1\lambda_2$ for some $\psi=(\zeta^i,\xi^j)\in\O_4$; the case of $\psi_-$ is analogous.
	We need to show that $Z=\{(\zeta^p,\xi^q)| (p,q)\in I\}$ for
	$I=\{ (i,j), (i + N, j),  (i, j + M),  (i + N, j + M),
	(2N - i, 2M - j), (N - i, 2M - j),  (2N - i, M - j),  (N - i, M - j)\}$.
	
	Now $\phi=(\zeta^p,\xi^q)\in Z$ if and only if $\phi_{+}=\psi_+$ or $\phi_{-}=\psi_+$, by~\eqref{eqn:lambdas}.
	Set $f(z,x)=\frac{1-z}{1-zx}$, so that $\phi_+=f(\zeta^{2p},\xi^{2q})$ and $\phi_-=f(\zeta^{-2p},\xi^{-2q})$.
	As these are roots of unity, we may compute in $\mathbb{C}$, where $f$ is injective on pairs $(z,x)$ of roots of unity with $z,x,zx\neq 1$: if $z=e^{\ima s}$ and $x=e^{\ima t}$, then $\frac{1}{f(z,x)}=e^{\ima t/2}\,\frac{\sin((s+t)/2)}{\sin(s/2)}$, whose argument determines $x$ and whose real factor then determines $z$.
	Hence $\phi_+=\psi_+$ if and only if $(\zeta^{2p},\xi^{2q})=(\zeta^{2i},\xi^{2j})$, and $\phi_-=\psi_+$ if and only if $(\zeta^{2p},\xi^{2q})=(\zeta^{-2i},\xi^{-2j})$.
	The index set $I$ then follows directly from the periodicities $N$ and $M$ of the squared generators, proving the claim and the proposition.
	\epf

	\subsection{Summary and graphical description}\label{sec:graphics}
	We have described simple modules of dimensions 1, 2 and 4; this has determined a partition:
	\begin{align}\label{eqn:Opart}
		\O=\O_1\sqcup \O_2\sqcup \O_4.
	\end{align}
	We now introduce a graphical perspective that allows a complete description of these modules as in~\eqref{eqn:simpledraw}.
	
	Recall that each one of these modules is generated by a component of a certain type $\chi\in\O$, of dimension 1. 
	We denote this component by $\ij$. 
	We write $\bij, \ibj$ and $\bibj$ for components of type $\overline{\chi}$, $-\overline{\chi}$ and $-\chi$, respectively.
	
	We use a horizontal labeled arrow $x\stackrel{a}{\rightarrow}y$ from vertex $x$ to vertex $y$ to represent the action of $a_1$, meaning $a_1\cdot x=a\,y$, $a\in\k$. When $a=1$, we omit the label. We use the same conventions for vertical arrows and the action of $a_2$, \emph{mutatis mutandis}. The absence of a horizontal (resp. vertical) arrow stands for the zero action of $a_1$ (resp. $a_2$) on that vertex. 
	We write $x\xleftrightarrow[b]{a}y$ to represent the assignments $a_1\cdot x=a\,y$ and $a_1\cdot y=b\,x$ (so $ab=\alpha_1$); similarly for vertical arrows and the action of $a_2$. Again, we omit the label when $a=1$ or $b=1$. 
	
	\begin{example}
		This notation allows us to define some indecomposable modules, as extensions of modules of dimensions 1 and 2. 
		\begin{enumerate}[leftmargin=*]
			\item[(a)] Fix $\chi\in\O_1$. We let $M_{1,h}^\chi\in\Ext^1(L_1^{\overline{\chi}},L_1^\chi)$ and $M_{1,v}^\chi\in\Ext^1(L_1^{-{\chi}},L_1^\chi)$ be the indecomposable modules:
			\begin{align}\label{eqn:ext1}
				M_{1,h}^\chi &: 
				\xymatrix{
					\ij&\bij\ar[l]}, 
				&
				M_{1,v}^\chi &: 
				\xymatrix{
					\ij\\
					\bibj.\ar[u]	}
			\end{align}
			For example, $M_{1,h}^\chi$ is the $H_\lambda$-module with basis $\{v,w\}$ such that $\lg v\rg_{|}\simeq S_\chi$ and $a_1\cdot w=v$ (also $a_1\cdot v=a_2\cdot v=a_2\cdot w=0$).
			\item[(b)] Fix $\chi\in\O_2$ and $a,b\in\k$, not both zero. When $\alpha_1\neq 0$ (respectively, $\alpha_2\neq 0$), we define $M_{2,h}^\chi(a,b)\in \Ext^1(L_{2,h}^{-{\chi}},L_{2,h}^\chi)$ (respectively, $M_{2,v}^\chi(a,b)\in \Ext^1(L_{2,v}^{\overline{\chi}},L_{2,v}^\chi)$) as follows:
			\begin{align}\label{eqn:ext2}
				M_{2,h}^\chi(a,b) &: 
				\xymatrix{
					\ij\ar@{<->}[r]_{\alpha_1} &\bij\\
					\bibj\ar@{<->}[r]_{\alpha_1}\ar[u]^{a} &\ibj\ar[u]^b
				}, &
				M_{2,v}^\chi(a,b) &: 
				\xymatrix{
					\ij\ar@{<->}[d]_{\alpha_2} &\bij\ar@{<->}[d]_{\alpha_2}\ar[l]^a\\
					\bibj &\ibj\ar[l]^b.}
			\end{align}
			When $a=0$ or $b=0$, we omit the corresponding arrow. We remark that $M_{2,-}^\chi(a,b)$ is the Baer sum $M_{2,-}^\chi(a,b)=aM_{2,-}^\chi(1,0)+bM_{2,-}^\chi(0,1)$.
		\end{enumerate}
	\end{example}
	
	\section{Simple modules and projective covers}\label{sec:classif}
	
	Fix the (projective) $H_\lambda$-module $P^\chi\coloneqq H_\lambda\ot_{\Gamma} S_\chi$. We shall combine the analysis of these modules with Lemma \ref{lem:quotients} to achieve the classification of simple modules. If $S_\chi=\lg z_{ij}\rg$, then we shall consider the induced basis
	$\{v_1,\dots,v_8\}$ given by 
	\begin{align}\begin{split}
			\label{eqn:basis-Pij}
			\{1\ot z_{ij},a_1\ot z_{ij},&a_2\ot z_{ij},a_1a_2\ot z_{ij},a_2a_1\ot z_{ij},\\
			&a_1a_2a_1\ot z_{ij},a_2a_1a_2\ot z_{ij},a_1a_2a_1a_2\ot z_{ij}\}.
		\end{split}
	\end{align}
	
	A straightforward computation leads to the following.
	\begin{lemma}
		As a $\Gamma$-module, $P^\chi=(S_{\chi})^2\oplus (S_{\overline{\chi}})^2\oplus (S_{-\overline{\chi}})^2\oplus (S_{-\chi})^2$. More precisely, in the basis~\eqref{eqn:basis-Pij}:
		\begin{align}\label{eqn:Pdecomposition}
			(P^\chi)_{|}&=S_{\chi}\oplus S_{\overline{\chi}}\oplus S_{-\chi}\oplus S_{-\overline{\chi}}\oplus S_{-\overline{\chi}}\oplus S_{-\chi}\oplus S_{\overline{\chi}}\oplus S_{\chi},
		\end{align}
		namely, the action of $\Gamma$ on this basis is given by the diagonal matrices:
		\begin{align*}
			[g_1]&=\zeta^{i}\diag(1, -1, -1, 1, 1, -1, -1, 1), \\
			[g_2]&=\xi^{j}\diag(1, 1, -1, -1, -1, -1, 1, 1).
		\end{align*}
		
		In turn, the action of $a_1$ and $a_2$ is given by the matrices:
		\begin{align}\label{eqn:matrices}
			[a_1]=\begin{bsmallmatrix*}
				0&\alpha_1&0&0&0&0&0&0\\
				1&0&0&0&0&0&0&0\\
				0&0&0&\alpha_1&0&0&0&0\\
				0&0&1&0&0&0&0&0\\
				0&0&0&0&0&\alpha_1&0&0\\
				0&0&0&0&1&0&0&0\\
				0&0&0&0&0&0&0&\alpha_1\\
				0&0&0&0&0&0&1&0
			\end{bsmallmatrix*},&&	
			[a_2]=\begin{bsmallmatrix*}
				0&0&\alpha_2&0&0&\alpha_3&0&0\\
				0&0&0&0&\alpha_2&0&0&0\\
				1&0&0&0&0&0&0&\alpha_3\\
				0&0&0&0&0&0&\alpha_2&0\\
				0&1&0&0&0&0&0&0\\
				0&0&0&0&0&0&0&-\alpha_2\\
				0&0&0&1&0&0&0&0\\
				0&0&0&0&0&-1&0&0
			\end{bsmallmatrix*}.
		\end{align}
	\end{lemma}
	
	\subsection{Four-dimensional submodules of $P^\chi$}\label{subsec:analysisP4}
	
	We investigate the possibility of having an irreducible submodule $L=L_4^\phi(d)\subset P ^\chi$, of dimension 4. Here $\phi\in\Omega(\chi)$. 
	Recall the discriminant $D=D(\chi)$ from equation~\eqref{eqn:discriminant}.
	
	\begin{lemma}\label{lem:proj-decomp}
		Let $\chi\in\O_4$ and set $P=P^\chi$. Define $\theta_{\pm}=\theta_{\pm}(\chi)$ by
		\begin{equation}\label{eqn:theta}
			\theta_{\pm}=\tfrac{-\alpha_3\pm\sqrt{D}}{2}.
		\end{equation}
		\begin{enumerate}[leftmargin=*]
			\item[(a)] If $D(\chi)\neq 0$ and $\alpha_1\alpha_2\neq 0$, then
			$P^\chi\simeq L_4^\chi(-\alpha_2/\theta_+)\oplus L_4^\chi(-\alpha_2/\theta_-)$.
			\item[(b)] If $D(\chi)\neq 0$ and $\alpha_1\alpha_2= 0$, then $P^\chi\simeq L_4^\chi(\alpha_2/\alpha_3)\oplus L_4^{\overline{\chi}}(\alpha_2/\alpha_3)$.
			\item[(c)] If $D(\chi)=0$, then there is a non-split extension
			\begin{align*}
				0 \to L_4^\chi(2\alpha_2/\alpha_3)\longhookrightarrow P^\chi \longtwoheadrightarrow L_4^\chi(2\alpha_2/\alpha_3)\to 0.
			\end{align*}
		\end{enumerate}
	\end{lemma}
	\pf
	We search for conditions under which a vector $w_1=\theta v_1+\gamma v_8$ generates a 4-dimensional submodule. As $\lg v_1\rg=P$, we can assume $\gamma=1$. 
	Hence $L$ should have the following basis:
	\begin{align}\label{eqn:basis-dim4}
		\begin{split}
			w_1&\coloneqq\theta v_1+v_8,\qquad \qquad \qquad \qquad 
			w_2\coloneqq a_1\cdot w_1=\theta v_2+\alpha_1 v_7,  \\
			w_3&\coloneqq a_2a_1\cdot w_1=\theta v_5+\alpha_1\alpha_2 v_4, \quad 
			w_4\coloneqq a_1a_2a_1\cdot w_1=\theta v_6+\alpha_1^2\alpha_2 v_3. 		
		\end{split}
	\end{align}
	We need $a_2\cdot w_4\in\k\{w_1\}$, and
	$
	a_2\cdot w_4=a_2a_1a_2a_1\cdot w_1=(\theta\alpha_3+\alpha_1^2\alpha_2^2)v_1-\theta v_8
	$.
	
	Observe that $\theta\neq 0$. Otherwise $w_1=v_8$, and either $\alpha_1\alpha_2\neq 0$, so that $a_2a_1a_2a_1\cdot w_1=\alpha_1^2\alpha_2^2 v_1\notin \k\{v_8\}$, or $\alpha_1\alpha_2=0$, so that $w_3=0$.
	Then
	\[
	a_2\cdot w_4=-\theta \left(\tfrac{\theta\alpha_3+\alpha_1^2\alpha_2^2}{-\theta}v_1+ v_8\right).
	\]
	Thus we need $\theta^2+\alpha_3\theta+\alpha_1^2\alpha_2^2=0$, that is, $\theta=\theta_{\pm}$ for $\theta_{\pm}$ as in~\eqref{eqn:theta}.
	
	In particular, $c=-\theta$ and it is easy to check that this defines a submodule $L_4^\chi(d)$, with $d=\dfrac{-\alpha_2}{\theta}$. Indeed, one checks $a_2\cdot w_1=-\frac{\alpha_2}{\theta}w_4$.
	
	(a) Whenever $\theta_+\neq \theta_-$ (namely, $D\neq 0$) and both are nonzero, this defines two linearly independent solutions $w_1^+=\theta_+ v_1+v_8$ and $w_1^-=\theta_- v_1+v_8$. The claim follows.
	
	(b) Assume $D\neq 0$. There is a single nonzero solution $\theta$ when $\alpha_1\alpha_2=0$. Hence $\theta=-\alpha_3$ and $L_4^\chi(\alpha_2/\alpha_3)\subset P^\chi$.
	Assume $\alpha_1=0$ (and thus $\alpha_3\neq 0$). If we set 
	\begin{align*}
		w_1&=v_7, & w_2&=a_1\cdot v_7=v_8, \\
		w_3&=a_2\cdot w_2=\alpha_3v_3-\alpha_2v_6, & w_4&=a_1\cdot w_3=\alpha_3v_4, 
	\end{align*}
	then we see that $a_2\cdot w_1=a_2\cdot v_7=\alpha_2v_4=\frac{\alpha_2}{\alpha_3}w_4$, which defines a submodule $\simeq L_4^{\overline{\chi}}(\alpha_2/\alpha_3)$.
	This proves the assertion in this case.
	Assume, alternatively, that $\alpha_1\neq 0$ and $\alpha_2=0$. Then $w_1=v_7$ again generates a submodule $\simeq L_4^{\overline{\chi}}(0)$ (observe that $a_2\cdot w_1=0$ in this case) and 
	$
	P^\chi\simeq L_4^\chi(0)\oplus L_4^{\overline{\chi}}(0)$; hence the claim holds.
	
	(c) When $D=0$, there is a single solution $\theta=-\alpha_3/2$ and we get that 
	$P^\chi/L_4^\chi(2\alpha_2/\alpha_3)\simeq L_4^\chi(2\alpha_2/\alpha_3)$. 
	Indeed, it is easy to check, using $\theta\bar{v_6}+\alpha_1^2\alpha_2\bar{v_3}=0$ in the quotient, where $\theta=-\alpha_3/2$, that $\lg \bar{v_1}\rg\simeq L_4^\chi(2\alpha_2/\alpha_3)$. Thus we obtain the extension from the statement.
	\epf
	
	\begin{remark}\label{rem:dim4-projectivematrices}
		Assume $D(\chi)=0$, let $w_1,w_2,w_3,w_4$ be as in~\eqref{eqn:basis-dim4} and consider the basis $\{w_1,w_2,w_3,w_4,v_1,v_2,v_5,v_6\}$ of $P^\chi$. It is an easy exercise to check that the actions of $a_1$ and $a_2$ are determined by the block matrices $[a_1]=\begin{psmallmatrix}
			A&0\\0&A
		\end{psmallmatrix}$ and $[a_2]=\begin{psmallmatrix}
			B&C\\0&B
		\end{psmallmatrix}$, 
		for $A$ and $B$ as in~\eqref{eqn:a1a2-dim4} and $C=\begin{psmallmatrix}
			0&0&0&-1\\0&0&0&0\\0&0&0&0\\ 1/(\alpha_1^2\alpha_2)&0&0&0
		\end{psmallmatrix}$.
	\end{remark}
	
	\subsection{The shape of $P^\chi$}\label{sec:cases}
	Fix $\chi\in\O$ and let $\alpha=(\alpha_1, \alpha_2, \alpha_3)\in\k^{3}$ be as in~\eqref{eqn:alpha}. In this part, we study the shape of $P^\chi$ according to how
	many of the $\alpha_i$ are zero. Then we have the following cases:
	\begin{align*}
		(i)\,& \alpha=(0, 0, 0), & (v)\,& \alpha=(0, \alpha_2, \alpha_3), \ \alpha_2\alpha_3\neq 0,\\
		(ii)\,& \alpha=(0, 0, \alpha_3), \ \alpha_3\neq 0, & (vi)\,& \alpha=(\alpha_1, 0, \alpha_3),  \ \alpha_1\alpha_3\neq 0,\\
		(iii)\,& \alpha=(0, \alpha_2, 0), \ \alpha_2\neq 0, & (vii)\,& \alpha=(\alpha_1, \alpha_2, 0), \ \alpha_1\alpha_2\neq 0,\\
		(iv)\,& \alpha=(\alpha_1, 0, 0), \ \alpha_1\neq 0, & (viii)\,& \alpha=(\alpha_1, \alpha_2, \alpha_3), \ \alpha_1\alpha_2\alpha_3\neq 0.
	\end{align*}
	
	\begin{proposition}\label{pro:projective}
		In each of the cases (i)--(viii), we have:
		\begin{enumerate}[leftmargin=*]
			\item[(i)] $P^\chi$ is indecomposable and the Jordan--H\"older series of $P^\chi=\lg v_1\rg$ is 
			\begin{align*}
				0\subset \lg v_8\rg\subset \lg v_7,-\rg\subset \lg v_6,-\rg\subset \lg v_5,-\rg\subset \lg v_4,-\rg\subset \lg v_3,-\rg\subset \lg v_2,-\rg\subset P^\chi
			\end{align*}
			with composition factors 
			$L_1^{\chi}, L_1^{\overline{\chi}}, L_1^{-\chi}, L_1^{-\overline{\chi}}, L_1^{-\overline{\chi}}, L_1^{-\chi}, L_1^{\overline{\chi}}, L_1^{\chi}$.
			\item[(ii)] $P^\chi\simeq L_{4}^{\chi}(0)\oplus L_{4}^{\overline{\chi}}(0)$.
			\item[(iii)]  $P^\chi$ is indecomposable and the Jordan--H\"older series of $P^\chi$ is 
			\begin{align*}
				0\subset \lg v_6, -v_8\rg\subset \lg -,v_2,v_5\rg\subset \lg -,v_4,v_7\rg\subset \lg -,v_1,v_3\rg= P^\chi
			\end{align*}
			with composition factors $L_{2,v}^{\chi}, L_{2,v}^{\overline{\chi}}, L_{2,v}^{\overline{\chi}}, L_{2,v}^{\chi}$.
			\item[(iv)] $P^\chi$ is indecomposable and the Jordan--H\"older series of $P^\chi$ is 
			\begin{align*}
				0\subset \lg v_7, v_8\rg\subset \lg -,v_3,v_4\rg\subset \lg -,v_5,v_6\rg\subset  \lg -,v_1,v_2\rg= P^\chi
			\end{align*}
			with composition factors $L_{2,h}^{\chi}, L_{2,h}^{-\chi}, L_{2,h}^{-\chi}, L_{2,h}^{\chi}$.
			\item[(v)] $P^\chi\simeq L_{4}^{\chi}(\alpha_2/\alpha_3)\oplus L_{4}^{\overline{\chi}}(\alpha_2/\alpha_3)$.
			\item[(vi)] $P^\chi\simeq L_{4}^{\chi}(0)\oplus L_{4}^{\overline{\chi}}(0)$.
			\item[(vii)] $P^\chi\simeq L_{4}^{\chi}(\sqrt{-1}/\alpha_1)\oplus L_{4}^{\chi}(-\sqrt{-1}/\alpha_1)$.
			\item[(viii)] If $D(\chi)\neq 0$, then $P^\chi\simeq L_{4}^{\chi}(-\alpha_2/\theta_+)\oplus L_{4}^{\chi}(-\alpha_2/\theta_-)$. When $D(\chi)=0$,  $P^\chi$ is indecomposable and the Jordan--H\"older series of $P^\chi$ is 
			\[
			0\subset L_4^\chi(2\alpha_2/\alpha_3)=\lg -\tfrac{\alpha_3}{2}v_1+v_8\rg \subset P^\chi
			\]
			with composition factors  $L_4^\chi(2\alpha_2/\alpha_3)$, twice.
		\end{enumerate}
	\end{proposition}
	\pf
	Items (ii) and (v)--(viii) follow from \S\ref{subsec:analysisP4}. For (viii), it remains to check that $P\coloneqq P^\chi$ is indecomposable. Suppose there is a decomposition $P=M\oplus N$; we can assume that $L\coloneqq L_4^\chi\subset M$. Now $L\simeq P/L\simeq M/L\oplus N$ implies that $N=0$ and $M=P$. 
	
	For (i), we also show that $P$ is indecomposable. Assume $P=M\oplus N$, $M\neq \{0\}$; then we can assume that there are $\alpha, \beta\in\k$ with $w=\alpha v_1+\beta v_8\in M$: namely the $\chi$-component of $M_{|}$ is non-trivial. We may assume $\beta\neq 0$, as otherwise $v_1\in M$ and $M=P$. If $\alpha\neq 0$, then $v_8=\alpha^{-1}a_1a_2a_1a_2\cdot w$, since $a_2\cdot v_8=0$. Hence $v_1\in M$ and $M=P$. The composition series follows by looking at ~\eqref{eqn:matrices} in this case; the same holds for chains in (iii) and (iv). 
	Indecomposability in cases (iii) and (iv) follows as in (i).
	\epf

	\subsection{Classification of simple modules}
	We present a complete classification of the irreducible representations of $H_\lambda$, for each $\lambda$. 
	We need to introduce some notation. 
	We consider the subsets $\H_N=\{\zeta^i:i\in\I_{N-1}\}\subset\G_{2N}$, $\H_M=\{\xi^j:j\in\I_{M-1}\}\subset\G_{2M}$. 
	We also define	
	\begin{align}\label{eqn:set}
		\s_{N,M}\coloneqq \{\chi=(\zeta^i,\xi^j)\in\O : \zeta^{2i}\xi^{2j}=1\}.
	\end{align}
	Observe that  $|\s_{N,M}|=4(N,M)$,	where $(N,M)$ stands for the greatest common divisor of $M$ and $N$. 
	
	\begin{theorem}\label{thm:simples}
		Fix $N,M\geq 1$, $\lambda=(\lambda_1,\lambda_2,\lambda_3)\in\k^3$. If  $L$ is an irreducible $H_\lambda$-module,
		then there exists $\chi \in \G_{2N} \times \G_{2M}$ such that $L$ is isomorphic to one of the following types: 
		$L\simeq L_1^\chi$, $L\simeq L_{2,h}^\chi$, $L\simeq L_{2,v}^\chi$ or $L\simeq L_4^\chi$. 
		More precisely, the simple $H_\lambda$-modules are, up to isomorphism, the following.
		\begin{enumerate}[leftmargin=*]
			\item If $\lambda_1=\lambda_2=\lambda_3=0$, then 
			\begin{itemize}[leftmargin=*]
				\item $4MN$ modules $L_1^\chi$, $\chi\in\O_1=\O$.
			\end{itemize}
			\item If $\lambda_1\neq 0$ and $\lambda_2=\lambda_3=0$, then
			\begin{itemize}[leftmargin=*]
				\item $4M$ modules $L_1^\chi$, $\chi\in \{\pm1\}\times\G_{2M}$.
				\item $2M(N-1)$ modules  $L_{2,h}^\chi$, $\chi\in\H_N\times\G_{2M}$. 
			\end{itemize}
			\item If $\lambda_2\neq 0$ and $\lambda_1=\lambda_3=0$, then
			\begin{itemize}[leftmargin=*]
				\item $4N$ modules  $L_1^\chi$, $\chi\in\G_{2N}\times \{\pm1\}$. 
				\item $2N(M-1)$ modules $L_{2,v}^\chi$, $\chi\in\G_{2N}\times \H_M$. 
			\end{itemize}
			\item If $\lambda_3\neq 0$ and $\lambda_1=\lambda_2=0$, then
			\begin{itemize}[leftmargin=*]
				\item $4(N,M)$ modules  $L_1^\chi$, $\chi\in\s_{N,M}$. 
				\item $4NM-4(N,M)$ modules $L_{4}^\chi$, $\chi\in\overline{\O}_4$.
			\end{itemize}
			\item If $\lambda_2\lambda_3\neq 0$ and $\lambda_1=0$, then
			\begin{itemize}[leftmargin=*]
				\item $4$ modules  $L_1^\chi$, $\chi\in\{\pm1\}^{\times 2}$.
				\item $2(N,M)-2$ modules  $L_{2,v}^\chi$, $\chi=(\chi_1,\chi_2)\in\s_{N,M}$, $\chi_2\in\H_M$. 
				\item $2NM-2(N,M)$ modules $L_4^\chi$, $\chi\in\overline{\O}_4$. 
			\end{itemize}
			\item If $\lambda_1\lambda_3\neq 0$ and $\lambda_2=0$, then
			\begin{itemize}[leftmargin=*]
				\item $4$ modules  $L_1^\chi$, $\chi\in\{\pm1\}^{\times 2}$.
				\item $2(N,M)-2$ modules  $L_{2,h}^\chi$, $\chi=(\chi_1,\chi_2)\in\s_{N,M}$, $\chi_1\in\H_N$.
				\item $2NM-2(N,M)$ modules $L_4^\chi$, $\chi\in\overline{\O}_4$. 
			\end{itemize}
			\item If $\lambda_1\lambda_2\neq 0$ and $\lambda_3=0$, then
			\begin{itemize}[leftmargin=*]
				\item 4 modules  $L_1^\chi$, $\chi\in\{\pm1\}^{\times 2}$.
				\item $2(M-1)$  modules  $L_{2,v}^\chi$, $\chi\in\{\pm1\}\times \H_M$. 
				\item $2M(N-1)$  modules $L_4^\chi$, $\chi\in\overline{\O}_4$.
			\end{itemize}
			\item If $\lambda_1\lambda_2\lambda_3\neq 0$, then
			\begin{itemize}[leftmargin=*]
				\item[(i)] when $\lambda_3=4\lambda_1\lambda_2$:
				\begin{itemize}[leftmargin=*]
					\item[$\bullet$] 4 modules  $L_1^\chi$, $\chi\in\{\pm1\}^{\times 2}$.
					\item[$\bullet$]  $2(N-1)$ modules $L_{2,h}^\chi$, $\chi\in\H_N\times \{\pm1\}$.
					\item[$\bullet$] $2N(M-1)$ modules $L_4^\chi$, $\chi\in\overline{\O}_4$.
				\end{itemize}
				\item[(ii)] when $\lambda_3=4\chi_\pm\lambda_1\lambda_2$ as in~\eqref{eqn:chipm}:
				\begin{itemize}[leftmargin=*]
					\item[$\bullet$] 4 modules  $L_1^\chi$, $\chi\in\{\pm1\}^{\times 2}$.
					\item[$\bullet$] $2(NM-2)$ modules $L_4^\chi$, $\chi\in\overline{\O}_4$.
				\end{itemize}
				\item[(iii)] when $\lambda_3\neq 4\lambda_1\lambda_2$, $\lambda_3\neq 4\chi_\pm\lambda_1\lambda_2$: 
				\begin{itemize}[leftmargin=*]
					\item[$\bullet$] 4 modules  $L_1^\chi$, $\chi\in\{\pm1\}^{\times 2}$.
					\item[$\bullet$] $2(NM-1)$ modules $L_4^\chi$, $\chi\in\overline{\O}_4$.
				\end{itemize}				
			\end{itemize}
		\end{enumerate}
	\end{theorem}
	\pf
	We apply Lemma \ref{lem:quotients}, by looking at the possible simple quotients of modules $P^\chi$, $\chi\in\O$, as described in Proposition \ref{pro:projective}. This shows that every simple module has dimension 1, 2, or 4, and all such modules
	have been described and classified above.
	
	The number of simple modules in each case follows by counting the subsets $\O_1$, $\O_2$, $\O_4\subset\O$. 
	We take into account the isomorphisms $L_{2,h}^\chi\simeq L_{2,h}^{\overline{\chi}}$ and $L_{2,v}^\chi\simeq L_{2,v}^{-\chi}$ to choose a representative in each case, so we have $\frac12|\O_2|$ modules of dimension 2.
	The same applies to the irreducible modules of dimension 4, following formula~\eqref{eqn:number=4} in Proposition \ref{pro:iso4}.
	
	We work out case \emph{(5)} as an example. In this case $\alpha_1=0$, $\alpha_2=\lambda_2(1-\xi^{2j})$ and $\alpha_3=\lambda_3(1-\zeta^{2i}\xi^{2j})$. Thus $\chi\in\O_1$ only if $\xi^{2j}=1$, so $\alpha_2=0$; hence $\alpha_3=\lambda_3(1-\zeta^{2i})$ and thus we also require $\zeta^{2i}=1$. 
	Hence $\O_1=\{\pm1\}^{\times 2}$. Now $\chi\in\O_2$ when $\alpha_2\neq0$, so $\xi^j\in\G_{2M}\setminus\{\pm1\}$ and $\alpha_3=0$, namely when $\chi\in\s_{N,M}$. 
	These are $|\s_{N,M}|-4$ possibilities which lead to $2(N,M)-2$ simple modules of dimension 2, up to isomorphism. The remaining $4NM-|\O_1|-|\O_2|=4NM-4-(4(N,M)-4)$ pairs $\chi\notin \O_1\cup\O_2$ give rise to the simple modules of dimension 4, of which there are $\frac{1}{2}(4NM-4(N,M))$ up to isomorphism.  
	\epf
	
	We finish this part with a connection to the theory of Hopf cocycles; recall the notation from \S\ref{sec:cocycles}. Keep in mind as well our \emph{convention} on p.~\pageref{page:convention} that $N,M\geq 2$ to avoid trivial deformations.
	
	\begin{corollary}\label{cor:cocycles}
		Let $\sigma$ be a Hopf cocycle such that $H_\lambda\simeq (\B_{\qb}\#\k\Gamma)_\sigma$.
		If there are simple $H_\lambda$-modules of dimensions 1, 2 and 4, then $\sigma$ is pure.
		
		Conversely, if $\sigma$ is pure, then there are simple $H_\lambda$-modules of dimensions 1, 2 and 4 if and only if $\lambda_3=0\neq \lambda_1\lambda_2$; or $\lambda_2 = 0 \neq \lambda_1\lambda_3$ and $(N, M) \neq 1$; or $\lambda_1 = 0 \neq \lambda_2\lambda_3$ and $(N, M)\neq 1$; or $\lambda_3 = 4\lambda_1\lambda_2 \neq 0$.
	\end{corollary}
	\pf
	According to Theorem \ref{thm:simples} above, there are simple $H_\lambda$-modules of dimensions 1, 2 and 4 only in cases \emph{(5)} to \emph{(8)}. In all of these cases, there are at least two nonzero parameters $\lambda_i$, $i=1,2,3$: hence the cocycle $\sigma$ is pure by~\cite[Theorem 6.13]{GS}.
	For the converse, we need to exclude the possibility $(N,M)=1$ in cases \emph{(5)} and \emph{(6)}, and $\lambda_3\neq 4\lambda_1\lambda_2$ in case \emph{(8)}.
	\epf

	\subsection{Projective covers} 
	By Proposition \ref{pro:projective}, we obtain the following.
	
	\begin{corollary}\label{cor:projectivecovers}
		We find the projective cover $P(L)$ of each simple module $L$.
		\begin{itemize}
			\item[(a)] The projective cover of $L_1^\chi$ is $P^\chi$.
			\item[(b)] The projective cover of $L_{2,\ast}^\chi$ is $P^\chi$. 
			\item[(c)] If $D(\chi)= 0$, then  $P(L_4^\chi(d))\simeq P^\chi$. Otherwise, $L_4^\chi(d)$ is projective.
		\end{itemize}
	\end{corollary}
	\pf
	(a) It follows from~\cite[Proposition 4.3 (ii)]{GI}, as $P^\chi$ is indecomposable in this case. 
	(b) We check the case $L_{2,h}^\chi$; the other case is analogous. Notice that we have $\alpha_2=\alpha_3=0$. 
	Recall also that we have a projection $P^\chi\twoheadrightarrow L_{2,h}^\chi$. Hence $P^\chi$ projects onto the projective cover $P(L_{2,h}^\chi)$ of $L_{2,h}^\chi$. 
	By Remark \ref{rem:dim-Hchi}, see~\eqref{eqn:dim32}, we have
	\[
	32=\dim L_{2,h}^\chi\dim P(L_{2,h}^\chi)+\dim L_{2,h}^{-\chi}\dim P(L_{2,h}^{-\chi})
	\]
	so $16=\dim P(L_{2,h}^\chi)+\dim P(L_{2,h}^{-\chi})$. As $L_{2,h}^{-\chi}\simeq L_{2,h}^{\chi}\ot L_1^{(-1,-1)}$, it follows that $\dim P(L_{2,h}^\chi)=\dim P(L_{2,h}^{-\chi})=8$ and thus $P(L_{2,h}^\chi)\simeq P^\chi$.
	
	(c) If $D(\chi)\neq 0$, then $L_4^\chi(d)$ is a direct summand of $P^{\chi'}$, for some $\chi'\in\Omega(\chi)$. 
	Indeed, if $\alpha_1\alpha_2=0$, then $d=\alpha_2/\alpha_3$ and $P^{\overline{\chi}}=L_4^{\overline{\chi}}(d)\oplus L_4^\chi(d)$, by Lemma \ref{lem:proj-decomp} (b). If $\alpha_1\alpha_2\neq 0$, then $d={(\alpha_3\pm\sqrt{D})}/{(2\alpha_1^2\alpha_2)}={-\theta_{\mp}}/{(\alpha_1^2\alpha_2)}$.
	Hence $L_4^\chi(d)=L_4^{\chi}(-\theta_{\mp}/(\alpha_1^2\alpha_2))\simeq L_4^{-\chi}(-\alpha_2/\theta_{\mp})$ by Proposition \ref{pro:iso4} (iii). Therefore, $L_4^\chi(d)$ is a direct summand of $P^{-\chi}$ by Lemma \ref{lem:proj-decomp} (a).
	
	Otherwise, there is a projection $P^\chi\twoheadrightarrow L_4^{\chi}$ and this map is essential as $L_4^\chi$ is the unique nonzero proper submodule of $P^\chi$ in this case.
	\epf

	\section{Extensions and the Gabriel quiver}\label{sec:gabriel}
	In this part we compute the extensions between simple modules, which allows us to determine that $H_\lambda$ is not of finite representation type.

	There are no nontrivial extensions between simple modules of different dimensions, by Lemma \ref{lem:indecomp}; hence $\dim\Ext^1(L,L')=0$ when $\dim L\neq \dim L'$. Also recall the interaction of $\Gamma$-isotypic components from~\eqref{eqn:squares}: this gives $\dim\Ext^1(L,L')=0$ when $L=L^\chi$ and $L'=L'^\phi$, with $\phi\notin \Omega(\chi)$. 
	
	We start with a lemma regarding extensions supported in $\O_4$ when $D(\chi)=0$. Recall from Lemma \ref{lem:proj-decomp} that the indecomposable projective module $P=P^\chi$ fits in an extension $0\to L_4\to P\to L_4\to 0$; see also Remark \ref{rem:dim4-projectivematrices}. 
	We  show that this nontrivial extension is unique, up to isomorphism.
	
	\begin{lemma}\label{lem:ext44}
		Let $\chi\in\O_4$ with $D(\chi)=0$ and set $L_4=L_4^\chi(2\alpha_2/\alpha_3)$. Let $r\geq 1$ and $M\not\simeq L_4^{\oplus\,r+1}$ be an $H_\lambda$-module fitting 
		in an extension $0\to L_4^{\oplus\,r}\to M\to L_4\to 0$. Then $M\simeq P^\chi\oplus N$, for some $H_\lambda$-module $N$.
		In particular, $\dim\Ext^1(L_4,L_4)=1$.	
	\end{lemma}
	\pf
	Let us fix an extension $0\to L_4^{\oplus r}\to M\to L_4\to 0$, $M\not\simeq L_4^{\oplus r+1}$.
	Fix a submodule $M'\subset M$ such that $M'\simeq L_4^{\oplus r}$. Let $u_1\in M_{|}[\chi]$ be such that the image $\bar{u}_1\in M/M'$ is not zero and $\lg \bar{u}_1\rg\simeq L_4$. We define $u_2\coloneqq a_1\cdot u_1$, $u_3\coloneqq a_2\cdot u_2$ and $u_4\coloneqq  a_1\cdot u_3$.
	
	Then there are $z\in  M'_{|}[-\chi]$ and $w\in  M'_{|}[\chi]$, not both zero, such that 
	\[
	a_2\cdot u_1=du_4+z, \quad a_2\cdot u_4=cu_1-w.
	\]
	Set $w_1\coloneqq w$, $w_2=a_1\cdot w_1$, $w_3=a_2\cdot w_2$, $w_4=a_1\cdot w_3$. As $w_1\in M'\simeq L_4^{\oplus r}$, it follows that $a_2\cdot w_1=dw_4$, namely $\lg w_1,w_2,w_3,w_4\rg\simeq L_4(d)$.
	
	On the other hand, as $dc=\alpha_2$, we have $a_2\cdot z=dw$, since
	\begin{align*}
		\alpha_2 u_1=a_2^2\cdot u_1=a_2\cdot(du_4+z)=dcu_1-dw+a_2\cdot z.
	\end{align*} 
	Namely, $z=\frac{d^2}{\alpha_2}w_4=(\frac{2\alpha_2}{\alpha_3})^2\frac{1}{\alpha_2}w_4=\frac{1}{\alpha_1^2\alpha_2}w_4$, since $\alpha_3^2=4\alpha_1^2\alpha_2^2$. 
	
	We check that $a_1a_2a_1a_2\cdot u_1=d\alpha_1^2\alpha_2 u_1+w_1$ and $a_2a_1a_2a_1\cdot u_1=cu_1-w_1$,
	so $(a_1a_2a_1a_2+a_2a_1a_2a_1)\cdot u_1=\alpha_3 u_1$. Hence we obtain a submodule $M''\subset M$ with basis $\{w_1,\dots,w_4,u_1,\dots,u_4\}$, which determines a non-split extension $0\to L_4\to M''\to L_4\to 0$ and for which the matrices $[a_1]$ and $[a_2]$ are as in Remark \ref{rem:dim4-projectivematrices}, so  $P^\chi\simeq M''\subset M$. As $P^\chi$ is injective (it is projective and $H_\lambda$ is a Hopf algebra), there is a submodule $N\subset M$ such that $M\simeq P^\chi\oplus N$.
	\epf
	
	\begin{proposition}\label{pro:extensions}
		Fix $\chi\in\O$ and let $\phi\in\Omega(\chi)$.
		The following holds.
		\begin{enumerate}[leftmargin=*]
			\item[(a)] If $\chi\in\O_1$, then $\dim\Ext^1(L_1^\phi,L_1 ^\chi)=\begin{cases}
				1, & \phi\in\{\overline{\chi},-\chi\};\\
				0, & \text{otherwise.}
			\end{cases}$
			\item[(b)] If $\chi\in\O_2$ and $\alpha_1\neq 0$, then $\dim\Ext^1(L_{2,h}^\phi,L_{2,h}^\chi)=\begin{cases}
				2, & \phi\in\{-\overline{\chi},-\chi\};\\
				0, & \text{otherwise.}
			\end{cases}$
			\item[(c)] If $\chi\in\O_2$ and $\alpha_2\neq 0$, then $\dim\Ext^1(L_{2,v}^\phi,L_{2,v}^\chi)=\begin{cases}
				2, & \phi\in\{-\overline{\chi},\overline{\chi}\};\\
				0, & \text{otherwise.}
			\end{cases}$
			\item[(d)] If $\chi\in\O_4$ and $D(\chi)\neq 0$, then $\dim\Ext^1(L_{4}^\phi(e),L_{4}^\chi(d))=0$. 
			\item[(e)] If $\chi\in\O_4$ and  $D(\chi)=0$, then 
			$\dim\Ext^1(L_{4}^\phi(e),L_{4}^\chi(d))=1$. 
		\end{enumerate}
	\end{proposition}
	
	\pf
	Case (d) follows since $L_4^\chi(d)$ is projective; case (e) is Lemma \ref{lem:ext44} for $r=1$.
	
	For case (a), let $M$ be an indecomposable module, necessarily of dimension 2, with $L_1^\chi\subset M$ and such that $M/L_1^\chi\simeq L_1^\phi$. This determines a basis $\{x,y\}$ with $\lg x\rg_{|}\simeq S_\chi$ and $\lg y\rg_{|}\simeq S_\phi$, together with $a_1\cdot x=a_2\cdot x=0$ and $a_1\cdot y,a_2\cdot y\in\k\{x\}$. As $\lg a_1\cdot y\rg_{|}\simeq S_{\overline{\phi}}$ and $\lg a_2\cdot y\rg_{|}\simeq S_{-{\phi}}$, we have that either $\phi=\overline{\chi}$ or $\phi=-\chi$. Thus $M$ is one of the modules in~\eqref{eqn:ext1}. 
	
	A similar analysis gives cases (b) and (c), for the modules in~\eqref{eqn:ext2}.
	\epf

	The following result follows from~\cite[Theorem X.2.6]{AuRS}, see also~\cite[Lemma 2.1]{GR}. We recall that the combination of these results states that an algebra whose separated quiver is not of finite Dynkin type cannot have finite representation type.

	\begin{corollary}
		$H_\lambda$ is not of finite representation type.
	\end{corollary}
	\pf
	By Proposition \ref{pro:extensions}, the Gabriel quiver of $H_\lambda$ is a disjoint union of quivers $G=G_1\sqcup G_2\sqcup G_4$, since there are no extensions between irreducible modules of different dimensions. 
	
	Now set $\widehat{\O}_k\coloneqq\widehat{\O}\cap\O_k$, $k\in\{1,2,4\}$; this is a set of representatives of the classes $\Omega(\chi)\subseteq\O_k$. Since there are no extensions between simple modules supported on different classes, $G_1=\bigsqcup_{\chi\in\widehat{\O}_1} G_1^\chi$ and $G_2=\bigsqcup_{\chi\in\widehat{\O}_2} G_2^\chi$ for
	\begin{align*}
		G_1^\chi: & \xymatrix{
			L_1^\chi\ar@/^.5pc/[r]\ar@/^.5pc/[d]&L_1^{\overline{\chi}}\ar@/^.5pc/[l]\ar@/^.5pc/[d]\\
			L_1^{-{\chi}}\ar@/^.5pc/[u]\ar@/^.5pc/[r] &L_1^{-\overline{\chi}}\ar@/^.5pc/[l]\ar@/^.5pc/[u]
		}
		&
		G_2^\chi:&
		\begin{cases}
			\xymatrix{L_{2,h}^\chi\ar@/^1pc/[r]\ar@/^.5pc/[r]&L_{2,h}^{-{\chi}}\ar@/^1pc/[l]\ar@/^.5pc/[l]}, \, \alpha_1\neq 0\\
			&\\
			\xymatrix{L_{2,v}^\chi\ar@/^1pc/[r]\ar@/^.5pc/[r]&L_{2,v}^{\overline{\chi}}\ar@/^1pc/[l]\ar@/^.5pc/[l]}, \, \alpha_2\neq 0.
		\end{cases} 
	\end{align*}
	Analogously,  $G_4=\bigsqcup_{\chi\in\widehat{\O}_4} G_4^\chi$ for
	\begin{align*}
		G_4^\chi:\begin{cases}
			\xymatrix{ L_4^\chi(\alpha_2/\alpha_3) \ L_4^{-\chi}(\alpha_2/\alpha_3)}, & \alpha_1\alpha_2=0, \\
			L_4^\chi(2\alpha_2/\alpha_3)\hspace*{-0.1cm}\xymatrix@1{ {\,}^{}\ar@(ur,dr)} , &  D(\chi)= 0,\\
			\xymatrix{ L_4^\chi(d) \ L_4^\chi(e)}, & \text{otherwise.}
		\end{cases}
	\end{align*}
	Thus, the separated quiver is the disjoint union $S=S_1\sqcup S_2\sqcup S_4$, for 
	\begin{align*}
		S_1&=\bigsqcup_{\mathclap{\chi\in\widehat{\O}_1}} (A_3^{(1)})^{\times 2}, & S_2&= \bigsqcup_{\mathclap{\chi\in\widehat{\O}_2}} (A_1^{(1)})^{\times 2}, & S_4&= \bigsqcup_{\mathclap{\substack{\chi\in\widehat{\O}_4, \\ D(\chi)\neq 0}}} A_1^{\times 4}
		\sqcup \bigsqcup_{\mathclap{\substack{\chi\in\widehat{\O}_4, \\ D(\chi)= 0}}} A_2.
	\end{align*}
	As this is not a quiver of finite type, the corollary follows. Here $A_n$, $A_n^{(1)}$, $n\in\N$, follow the standard notation for Dynkin and Euclidean diagrams and $D^{\times k}$ stands for the diagram $\underbrace{D\sqcup\dots\sqcup D}_{k\text{ times}}$. 
	\epf
	
	\section{Indecomposable modules}\label{sec:indecomposables}
	
	As established previously, the algebras $H_\lambda$ are not of finite representation type. In this section, we investigate the structure of their indecomposable modules. Any such module must be supported on a single component $\O_k$, for some $k\in\{1,2,4\}$. We classify indecomposable modules supported on $\O_4$ in \S\ref{sec:indec4}. For modules supported on $\O_1\cup\O_2$, we classify indecomposable modules of small dimension in \S\ref{sec:small-dim}; these will be useful for computations in Section \ref{sec:quotient}. Finally, in \S\ref{sec:any-dim}, we construct an indecomposable $H_\lambda$-module of dimension $n$ for each $n\in\N$. Section \ref{sec:cut-off} includes a criterion for identifying projective modules supported on $\O_1\cup\O_2$.
	
	\subsection{Indecomposable modules supported on $\O_4$}\label{sec:indec4}
	
	We present a classification of indecomposable modules supported on $\O_4$.
	\begin{proposition}\label{pro:indec4}
		Let $M=M^\chi$ be an indecomposable $H_\lambda$-module, 	$\chi\in\O_4$. 
		\begin{enumerate}
			\item[(a)] If $D(\chi)\neq 0$, then $M\simeq L_4^\chi$. 
			\item[(b)] If $D(\chi)=0$, then $M\simeq L_4^\chi$ or $M\simeq P^\chi$.
		\end{enumerate}
	\end{proposition}
	\pf
	$(a)$ follows, as every simple submodule $L\subset M$ is projective.
	As for $(b)$, assume $M$ is not simple. Then the Jordan--H\"older series of $M$ necessarily contains a segment of the form $0\to L_4\to\dots \to M_1\simeq L_4^{\oplus\,r}\to M_2\subseteq M$, $r\geq 1$, with $M_2/M_1\simeq L_4$ and $M_2\not\simeq L_4^{\oplus\,r+1}$. By Lemma \ref{lem:ext44}, $M_2$ splits as $M_2\simeq P^\chi\oplus L_4^{\oplus\,r-1}$. A similar decomposition thus holds for $M$, namely $M\simeq P^\chi\oplus M'$ because $P^\chi$ is injective: a contradiction unless $M\simeq P^\chi$.
	\epf
	
	\subsection{Cut-off lemmas}\label{sec:cut-off}
	The following two lemmas provide a useful characterization of indecomposable modules $M=M^\chi$, with $\chi\in\O_i$, $i\in\I_2$.
	
	\begin{lemma}\label{lem:arrow41}
		Fix $\chi\in\O_1$ and let $M=M^\chi$ be an indecomposable $H_\lambda$-module. 
		If there is $x\in M$ such that $a_1a_2a_1a_2\cdot x\neq 0$, then $M\simeq P^\chi$.
	\end{lemma}
	\pf
	Recall that we have a $\Gamma$-decomposition $M_{|}=M_{|}[\chi]\oplus M_{|}[\overline{\chi}]\oplus M_{|}[-\overline{\chi}]\oplus M_{|}[-\chi]$. We can assume, without loss of generality, that $x\in M_{|}[\chi]$. 
	
	Now, as $-a_2a_1a_2a_1\cdot x=a_1a_2a_1a_2\cdot x\neq 0$, we have nonzero elements:
	\begin{align*}
		x_1&\coloneqq x, &y_1&\coloneqq a_1\cdot x, & w_1&\coloneqq a_2\cdot x, &z_1&\coloneqq a_1a_2\cdot x,\\
		x_2&\coloneqq a_1a_2a_1a_2\cdot x, &y_2&\coloneqq a_2a_1a_2 \cdot x, & w_2&\coloneqq a_1a_2a_1\cdot  x, &z_2&\coloneqq a_2a_1\cdot x.
	\end{align*}
	Observe that the subsets $\{x_1,x_2\}$, $\{y_1,y_2\}$, $\{w_1,w_2\}$, $\{z_1,z_2\}$ belong to different isotypic components of $M_{|}$. Moreover, each subset is linearly independent. Indeed, for the first one, we have $a_1\cdot x_2=0$ and $a_1\cdot x_1=y_1\neq 0$. The others follow similarly. This determines an 8-dimensional subspace $M'\subset M$ with basis $\{x_1,  y_1, w_1,z_1,z_2,  w_2, y_2,x_2\}$. 
	Moreover, this is a submodule, isomorphic to $P^\chi$. Indeed, for this (ordered) basis, $M'_{|}$ is as in~\eqref{eqn:Pdecomposition} and the matrices $[a_1]$ and $[a_2]$ coincide with those in~\eqref{eqn:matrices}, for the basis $\{v_1,\dots,v_8\}$ there. As $P^\chi$ is injective, this determines a complement $M''\subset M$ so that $M\simeq P^\chi\oplus M''$, which is a contradiction unless $M\simeq P^\chi$ as stated.
	\epf

	An analogous characterization holds when $M=M^\chi$, with $\chi\in\O_2$.
	
	\begin{lemma}\label{lem:arrow42}
		Fix $\chi\in\O_2$ and let $M=M^\chi$ be an indecomposable $H_\lambda$-module. 
		If there is $x\in M$ such that $a_1a_2a_1a_2\cdot x\neq 0$, then $M\simeq P^\chi$.
	\end{lemma}
	\pf
	Assume $\alpha_1\neq 0$, $\alpha_2=0$; the symmetric case is equivalent. Now, the result follows as in Lemma \ref{lem:arrow41}: once again we obtain four subsets $\{x_1,x_2\}$, $\dots$,  $\{z_1,z_2\}$ located in different $\Gamma$-components of $M$. Each subset is linearly independent as one of the elements is annihilated by $a_2$ and the other is not. We consider the subspace $M'\subset M$ generated  by $\{x_1,\dots,z_2\}$ and check that this defines a submodule for which the matrices $[a_1]$ and $[a_2]$ are as in~\eqref{eqn:matrices}. The proof ends as in Lemma \ref{lem:arrow41}.
	\epf
	
	\subsection{Indecomposable modules of small dimension}\label{sec:small-dim}
	We classify indecomposable modules of dimensions 3 and 4.
	
	\subsubsection{Indecomposable modules of dimension 3}\label{sec:indec3}
	
	We begin by listing all indecomposable modules of dimension 3. Notice that such a module $M$ necessarily satisfies $M=M^\chi$, with $\chi\in\O_1$, by Lemma \ref{lem:indecomp} and Example \ref{exa:qdim-simple}.
	
	\begin{definition}\label{def:indec3}
		For $\chi\in\O_1$, we introduce the following collection of (indecomposable) modules:
		\tiny{	\begin{align*}
				M^\chi_{3,1}:\xymatrix{
					\ij\ar@{<-}[r] &\bij\\
					&\ibj\ar[u]
				},  
				&& 
				M^\chi_{3,2}:\xymatrix{
					\ij\ar@{<-}[r] &\bij\ar[d]\\
					&\ibj
				},  
				&&
				M^\chi_{3,3}:\xymatrix{
					\ij\ar@{<-}[r] &\bij\\
					\bibj\ar[u] &
				},
				&&
				M^\chi_{3,4}:\xymatrix{
					\ij & \\
					\bibj\ar@{<-}[r]\ar[u] &\ibj
				}.
		\end{align*}}
	\end{definition}

	\begin{lemma}\label{lem:indecomp3}
		Let $M$ be an indecomposable $H_\lambda$-module of dimension 3. 
		Then there is a unique pair  $(\chi,i)\in \O_1\times \I_4$ such that $M\simeq M_{3,i}^\chi$. 
		
		Moreover, $(M_{3,1}^\chi)^\ast\simeq M_{3,4}^{-\overline{\chi^{-1}}}$ and $(M_{3,2}^\chi)^\ast\simeq M_{3,3}^{-\overline{\chi^{-1}}}$.	
	\end{lemma}
	\pf
	Let $M$ be such a module; we may assume that there is $x\in M$ with $L_1^\chi\simeq \lg x\rg\subset M$, and $\dim M=3$. In particular, note that $a_1\cdot x=a_2 \cdot x=0$.
	
	We look at $M_{|}= S_{\chi} \oplus S_{\chi'} \oplus S_{\chi''}$. Fix $y,z$ such that $\lg y \rg_{|}\simeq S_{\chi'}$ and $\lg z \rg_{|}\simeq S_{\chi''}$, and $\{x, y, z\}$ is linearly independent.
	
	\begin{claim}
		$\chi, \chi', \chi''$ are pairwise different.
	\end{claim}
	
	Assume first that $\chi'=\chi$. If $a_1\cdot y\neq 0$ and $a_2\cdot y\neq 0$, then $\{x, y, a_1\cdot y, a_2\cdot y\}$ is linearly independent, a contradiction.
	If $a_1\cdot y=a_2\cdot y=0$ then:
	\begin{itemize}[leftmargin=*]
		\item if $\chi''=\chi$ we obtain a decomposition for $M$ as $(L_1^\chi)^3$, a contradiction;
		\item if $\chi''=-\overline{\chi}$, then $M\simeq (L_1^\chi)^2\oplus L_1^{\chi''}$, a contradiction. The same holds for $\chi''=-\chi$.
		\item if $\chi''=\overline{\chi}$, then $a_2\cdot z=0$, and $a_1\cdot z \in \k\{x, y\}$. If $a_1\cdot z=0$, then $M\simeq (L_1^\chi)^2\oplus L_1^{\chi''}$. If $a_1\cdot z\neq 0$, a change of basis gives $M\simeq \lg z, a_1\cdot z\rg \oplus L_1^\chi$.
	\end{itemize}
	
	On the other hand, assume $\chi' = \chi''$:
	\begin{itemize}[leftmargin=*]
		\item if $\chi'=\chi''=-\overline{\chi}$ then we get $M\simeq L_1^{\chi} \oplus \lg y, z\rg$, a contradiction.
		\item if $\chi'=\chi''=\overline{\chi}$, by a dimension argument, $a_2\cdot y=a_2\cdot z=0$, and similarly for $\chi'=\chi''=-\chi$ when $a_2\cdot y=0=a_2\cdot z$ (if it does not occur, i.e., $a_2\cdot y\neq0$ or $a_2\cdot z\neq0$, then $M$ can be decomposed as a direct sum).
		Now,
		\begin{itemize}[leftmargin=*]
			\item if $a_1\cdot y=0$ or $a_1\cdot z=0$, then $M$ is decomposable;
			\item if $a_1 \cdot y=\lambda x$ and $a_1 \cdot z=\mu x$, $\lambda,\mu\in\k^{\times}$, then $M\simeq \lg x, y \rg \oplus \lg y-\frac{\lambda}{\mu}z\rg$, a contradiction.
		\end{itemize}
	\end{itemize}
	This shows the claim.
	
	We have the following possibilities for the (ordered) triple $(\chi, \chi', \chi'')$:
	\begin{align*}
		1)\, (\chi, \overline{\chi}, -\overline{\chi}), && 2)\, (\chi, -\chi, -\overline{\chi}), && 3)\, (\chi, \overline{\chi}, -\chi).
	\end{align*}
	
	For case $1)$, recall that $a_1\cdot x=a_2\cdot x=0$. If $a_1\cdot y= 0$, then we obtain a direct summand $L_1^\chi$ $(\simeq \lg x \rg)$ for $M$. Assume $a_1\cdot y=x$. Also, $a_1\cdot z=0$ necessarily (it would be of type $-\chi$).
	If $a_2\cdot y \neq 0$, we can assume, up to rescaling, $a_2\cdot y=z$ (which implies $a_2\cdot z=0$) and thus $M\simeq M_{3,2}^\chi$.
	If $a_2\cdot y=0$, then $a_2 \cdot z\neq 0$ (otherwise, $\lg z\rg$ is a direct summand of $M$). We can assume $a_2 \cdot z=y$. Thus, $M\simeq M_{3,1}^\chi$.
	
	Case $2)$ is similar and we obtain either $M\simeq M_{3,4}^\chi$ or $M\simeq M_{3,2}^{-\overline{\chi}}$.
	
	For $3)$, $a_2\cdot y=0$ (type $\overline{\chi}$) and $a_1\cdot z=0$ (type $-\chi$). If $a_1\cdot y=0$ or $a_2\cdot z=0$, we find a direct summand for $M$ (contradiction). Then, $a_1\cdot y=\lambda x$ and $a_2\cdot z=\mu x$, $\lambda,\mu\in\k^{\times}$. Setting $\bar{y}=\frac{1}{\lambda}y$ and $\bar{z}=\frac{1}{\mu}z$, we get $M\simeq M_{3,3}^{\chi}$.
	\epf

	\subsubsection{Indecomposable modules of dimension 4}\label{sec:indec4-class}
	
	In Section \ref{sec:quotient} we shall define a quotient category $\uRep H_\lambda$. To get a glimpse of its fusion rules, we shall compute the tensor products of the modules in Lemma \ref{lem:indecomp3}. We shall encounter some of the following modules of dimension 4 as direct summands.
	\begin{definition}\label{def:c4}
		For each $\chi\in\O_1$, we introduce the following (families of) indecomposable $H_\lambda$-modules of dimension 4; for $\mu\in\k^\times$:
		\tiny{\begin{align*}
				C_{\underline{3},\mu}^{\chi} :\hspace*{-.2cm} 
				\xymatrix{
					\ij &\bij\ar[l]\\
					\bibj\ar[u]^{\mu}\ar[r] &\ibj\ar[u]
				},
				&& 
				C_{2,\mu}^{\chi} : \hspace*{-.2cm}
				\xymatrix{
					\ij &\bij\ar[l]\\
					\bibj\ar[u]^{\mu} &\ibj\ar[u]\ar[l]
				}, && 
				C_{\bar{3},\mu}^{\chi} :\hspace*{-.2cm} 
				\xymatrix{
					\ij &\bij\ar[l]^{\mu}\ar[d]\\
					\bibj\ar[u] &\ibj\ar[l]
				},
				&&
				C_{1,\mu}^{\chi} : \hspace*{-.2cm}
				\xymatrix{
					\ij &\bij\ar[l]\ar[d]\\
					\bibj\ar[u]^{\mu}\ar[r] &\ibj}.
		\end{align*}}
		\normalsize
		These modules are pairwise non-isomorphic, except for $C_{1,\mu}^{\chi}\simeq C_{1,\mu}^{-\overline{\chi}}$.
	\end{definition}
	
	\begin{lemma}\label{lem:indecomp4}
		Let $M=M^\chi$ be an indecomposable module of dimension 4. 
		\begin{itemize}[leftmargin=*]
			\item If $\chi\in\O_4$, then $M\simeq L_4^\chi(d)$ is simple as in Proposition \ref{pro:dim4}.
			\item If $\chi\in\O_2$, then $M\simeq M_{2,h}^\chi(a,b)$ or $M\simeq M_{2,v}^\chi(a,b)$ as in~\eqref{eqn:ext2}.
			\item If  $\chi\in\O_1$, then either $M\simeq C_{\ast,\mu}^{\chi}$ for some $\mu\in\k^\times$, as in Definition \ref{def:c4}, or $M$ is isomorphic to one (and only one) of the following:
			\tiny{		\begin{align*}
					M^\chi_{4,1}:\hspace*{-.2cm}\xymatrix{
						\ij &\bij\ar[l]\\
						\bibj\ar[r] &\ibj\ar[u]
					},  
					&&
					M^\chi_{4,2}:\hspace*{-.2cm}\xymatrix{
						\ij &\bij\ar[d]\\
						\bibj\ar[u] & \ibj\ar[l]
					},  
					&&
					M^\chi_{4,3}:\hspace*{-.2cm}\xymatrix{
						\ij &\bij\ar[l]\ar[d]\\
						\bibj &\ibj\ar[l]
					},
					&&
					M^\chi_{4,4}:\hspace*{-.2cm}\xymatrix{
						\ij & \bij\ar[l]\\
						\bibj\ar[u] &\ibj\ar[u]
					},
					\\
					M^\chi_{4,5}:\hspace*{-.2cm}\xymatrix{
						\ij &\bij\\
						\bibj\ar[u] &\ibj\ar[u]\ar[l]
					},  
					&&
					M^\chi_{4,6}:\hspace*{-.2cm}\xymatrix{
						\ij &\bij\ar[l]\\
						\bibj\ar[u] & \ibj\ar[l]
					},  
					&&
					M^\chi_{4,7}:\hspace*{-.2cm}\xymatrix{
						\ij &\bij\ar[l]\ar[d]\\
						\bibj\ar[u] &\ibj
					},
					&&
					M^\chi_{4,8}:\hspace*{-.2cm}\xymatrix{
						\ij & \bij\ar[l]\\
						\bibj\ar[r]\ar[u] &\ibj
					}.		
			\end{align*}}
		\end{itemize}
	\end{lemma}
	\pf
	Analogous to Lemma \ref{lem:indecomp3}.
	\epf
	
	\subsection{Indecomposable modules of arbitrary dimension}\label{sec:any-dim}
	
	For each $\chi\in\O_1$ and each $n\in\N$, and for each $\chi\in\O_2$ and each even $n\in\N$, we construct an indecomposable $H_\lambda$-module $M=M^\chi$ of dimension $n$.
	
	\begin{definition}\label{def:indecn-1}
		Fix $\chi\in\O_1$ and $n,k\in\N$ such that $n=4k$. 
		Define $Q_n^\chi$ as the $H_\lambda$-module with basis 
		$
		\{x_{1,1},x_{2,1},x_{3,1},x_{4,1},x_{1,2},x_{2,2},\dots, x_{1,k},x_{2,k},x_{3,k},x_{4,k}\}
		$
		so that $x_{1,j}, x_{2,j}, x_{3,j}, x_{4,j}$, $j\in\I_k$, have weights $\chi, \overline{\chi}, -\overline{\chi}$, $-\chi$, respectively, and the action is given as follows, for $j\in\I_k$, with the convention $x_{1,k+1}\coloneqq0$:
		\begin{align*}
			a_1\cdot x_{1,j}&=0, & a_2\cdot x_{1,j}&=0, & a_1\cdot x_{2,j}&=x_{1,j}, & a_2\cdot x_{2,j}&=0,\\
			a_1\cdot x_{3,j}&=0, & a_2\cdot x_{3,j}&=x_{2,j}, & a_1\cdot x_{4,j}&=x_{3,j}, & a_2\cdot x_{4,j}&=x_{1,j+1}.
		\end{align*}
		In turn, if $n=4k+r$ with $1\leq r\leq 3$, we set $Q_n^\chi$ to be the submodule of $Q^\chi_{4(k+1)}$ generated by the first $n$ elements in the ordered basis.
	\end{definition}
	
	\begin{example}
		We have already encountered some of these modules, since
		$Q_1^\chi=L_1^\chi$, $Q_2^\chi=M_{1,h}^\chi$ as in~\eqref{eqn:ext1}, $Q_3^\chi=M^\chi_{3,1}$ as in Definition \ref{def:indec3}
		and  $Q_4^\chi=M_{4,1}^\chi$  as in Lemma \ref{lem:indecomp4}. 
	\end{example}
	
	See Remark \ref{rem:sketch} for a visual description.
	
	\begin{lemma}\label{lem:indecn-1}
		For each $n\in\N$, $\chi\in\O_1$, $Q_n^\chi$ is an indecomposable $H_\lambda$-module of dimension $n$. 
	\end{lemma}
	\pf
	We deal with the case $n=4k$; the other cases are analogous.
	Set $M=Q_n^\chi$. We have $\dim M=n$ by definition; it is also clear that $a_1^2=a_2^2=a_1a_2a_1a_2=a_2a_1a_2a_1=0$ on $M$, which determines an $H_\lambda$-module structure.
	Moreover, the socle $\soc M$ of $M$ is the submodule $\lg x_{1,1},\dots, x_{1,k}\rg\simeq (L_1^\chi)^{\oplus\,k}$.
	Assume now that $M=N_1\oplus N_2$ and $x=x_{1,1}+\sum_{j=2}^{k}\lambda_jx_{1,j}\in N_1$, for some $\lambda_j\in\k$, $j>1$. 
	We claim that $y=x_{2,1}+\sum_{j=2}^{k}\lambda_jx_{2,j}\in N_1$. 
	Indeed, if $y=y_1+y_2$, $y_i\in N_i\cap\lg x_{2,1},\dots, x_{2,k} \rg$, $i\in\I_2$, then as $x=a_1\cdot y\in N_1$ we get that $a_1\cdot y_2=0$.
	This gives $L_1^{\overline{\chi}}\simeq \lg y_2\rg\subseteq M$, a contradiction. 
	This argument also shows that $z=x_{3,1}+\sum_{j=2}^{k}\lambda_jx_{3,j}\in N_1$ and $w=x_{4,1}+\sum_{j=2}^{k}\lambda_jx_{4,j}\in N_1$. Next, this implies that $N_1\ni a_2\cdot w=x_{1,2}+\sum_{j=2}^{k-1}\lambda_jx_{1,j+1}$. A recursive application of this procedure leads to $x_{1,k}\in N_1$. 
	
	On the one hand, if $y_1'\in N_1$ and $y_2'=\sum\eta_j x_{2,j}\in N_2$ are such that $a_1\cdot(y_1'+y_2')=x_{1,k}$, then $y_2'=0$ and $y_1'=x_{2,k}\in N_1$. Similarly, $x_{3,k}, x_{4,k}\in N_1$.
	
	On the other hand, if $w_1\in N_1$ and $w_2=\sum\mu_j x_{4,j}\in N_2$ are such that $a_2\cdot (w_1+w_2)=x_{1,k}$, we get that $\mu_j=0$, $j\in \I_{k-1}$; hence $w_2=\mu_k x_{4,k}\in N_2$ so $w_2=0$. This implies that $x_{4,k-1}\in N_1$, so $x_{3,k-1},x_{2,k-1}\in N_1$. Following this path, we get $x_{i,\ell}\in N_1$, $\ell\in\I_{k-1}$.
	
	Therefore $M=N_1$ and $N_2=\{0\}$; hence $M$ is indecomposable.
	\epf
	
	\begin{definition}\label{def:indecn-2}
		Fix $\chi\in\O_2$ with $\alpha_1\neq 0$, and let $n,k\in\N$ be such that $n=4k$. 
		We define $Q_{n,h}^\chi$ as the $H_\lambda$-module with basis given by the set 
		$
		\{x_{1,1},y_{1,1},x_{2,1},y_{2,1},\dots, x_{1,k},y_{1,k},x_{2,k},y_{2,k}\}
		$
		and action:
		\begin{align*}
			a_1\cdot x_{1,j}&=y_{1,j}, & a_2\cdot x_{1,j}&=0, & a_1\cdot y_{1,j}&=\alpha_1x_{1,j}, & a_2\cdot y_{1,j}&=0,\\
			a_1\cdot x_{2,j}&=y_{2,j}, & a_2\cdot x_{2,j}&=x_{1,j}, & a_1\cdot y_{2,j}&=\alpha_1 x_{2,j}, & a_2\cdot y_{2,j}&=y_{1,j+1}.
		\end{align*}
		When $n=4k+2$, we let  $Q_{n,h}^\chi$ be the submodule of $Q_{4(k+1),h}^\chi$ generated by the first $n$ basis elements $x_{1,1},y_{1,1},x_{2,1},y_{2,1},\dots, x_{1,k+1},y_{1,k+1}$.
		
		Analogously, when $\alpha_2\neq 0=\alpha_1$, we define $Q_{n,v}^\chi$ as the $H_\lambda$-module with linear basis 
		$
		\{x_{1,1},y_{1,1},x_{2,1},y_{2,1},\dots, x_{1,k},y_{1,k},x_{2,k},y_{2,k}\}
		$
		and action:
		\begin{align*}
			a_1\cdot x_{1,j}&=0,       & a_2\cdot x_{1,j}&=y_{1,j}, & a_1\cdot y_{1,j}&= 0,  & a_2\cdot y_{1,j}&=\alpha_2 x_{1,j},\\
			a_1\cdot x_{2,j}&=x_{1,j}, & a_2\cdot x_{2,j}&=y_{2,j}, & a_1\cdot y_{2,j}&= y_{1,j+1}, & a_2\cdot y_{2,j}&=\alpha_2 x_{2,j}.
		\end{align*}
		When $n=4k+2$, we let  $Q_{n,v}^\chi$ be the submodule of $Q_{4(k+1),v}^\chi$ generated by the first $n$ basis elements $x_{1,1},y_{1,1},x_{2,1},y_{2,1},\dots, x_{1,k+1},y_{1,k+1}$.
	\end{definition}
	
	\begin{example}
		$Q_{2,\ast}^\chi=L_{2,\ast}^\chi$ and $Q_{4,\ast}^\chi=M_{2,\ast}^\chi(1,0)$ as in ~\eqref{eqn:ext2}.
	\end{example}
	
	\begin{lemma}\label{lem:indecn-2}
		For each even $n\in\N$, $\chi\in\O_2$, $Q_{n,h}^\chi$ and $Q_{n,v}^\chi$ are indecomposable $H_\lambda$-modules of dimension $n$. 
	\end{lemma}
	\pf
	It follows as in Lemma \ref{lem:indecn-1}; here $\soc Q_{n,h}^\chi\simeq (L_{2,h}^\chi)^{\oplus k}$. 
	\epf
	
	\begin{remark}\label{rem:sketch}
		Let us sketch the modules $Q_{4k}^\chi$ and $Q_{4k,h}^\chi$. We add axes to stress the different components in~\eqref{eqn:squares}:
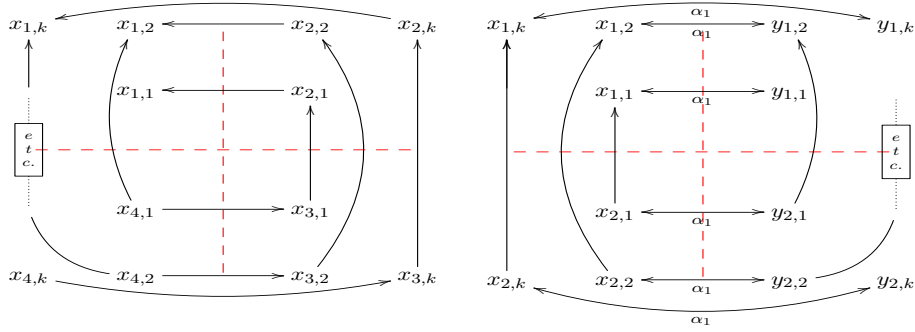
		\begin{align*}
			\resizebox{6cm}{.12cm}{
				\xymatrix{
					{x_{1,k}}&{x_{1,2}} &\ar@{--}@[red][dddd]&{x_{2,2}}\ar[ll]&{x_{2,k}}\ar@/_1pc/[llll]\\
					\ar[u]&{x_{1,1}} &&{x_{2,1}}\ar[ll]&\\
					\ar@{--}@[red][rrrr]&&&&\\
					\ar@{.}[uu]|*+[F-]
					{\tiny
						\begin{matrix}
							e\\t\\c.
					\end{matrix}}
					&{x_{4,1}}\ar@/^1pc/[uuu]\ar[rr]&&{x_{3,1}}\ar[uu]&\\
					{x_{4,k}}\ar@/_1pc/[rrrr]&{x_{4,2}}\ar[rr]\ar@/^1pc/@{-}[ul]&&{x_{3,2}}\ar@/_2pc/[uuuu]&{x_{3,k}}\ar[uuuu]
				}
			}
			&&
			\resizebox{6cm}{.28cm}{
				\xymatrix{
					{x_{1,k}}&{x_{1,2}} &\ar@{--}@[red][dddd]&{y_{1,2}}\ar@{<->}^{\alpha_1}[ll]&{y_{1,k}}\ar@{<->}^{\alpha_1}@/_1pc/[llll]\\
					\ar[u]&{x_{1,1}} &&{y_{1,1}}\ar@{<->}^{\alpha_1}[ll]&\\
					\ar@{--}@[red][rrrr]&&&&\\
					&{x_{2,1}}\ar[uu]\ar@{<->}_{\alpha_1}[rr]&&{y_{2,1}}\ar@/_1pc/[uuu]&\ar@{.}[uu]|*+[F-]
					{\tiny
						\begin{matrix}
							e\\t\\c.
					\end{matrix}}\\
					{x_{2,k}}\ar[uuuu]\ar@{<->}_{\alpha_1}@/_1.5pc/[rrrr]&{x_{2,2}}\ar@{<->}_{\alpha_1}[rr]\ar@/^2pc/[uuuu]&&{y_{2,2}}\ar@/_1pc/@{-}[ur]&{y_{2,k}}
				}
			}
		\end{align*}
		The ``vertical'' case $Q_{4k,v}^\chi$ is obtained by \emph{flipping} (rotating $90^\circ$ to the left and reflecting along the horizontal axis) the second diagram.
	\end{remark}

	\section{A spherical category}\label{sec:quotient}
	
	A spherical Hopf algebra~\cite{BaW-adv} is a pair $(H, \omega)$, where $H$ is a Hopf algebra and  $\omega\in G(H)$ is such that
	\begin{align}\label{eq:omega-square-antipode}
		S^2(x) &= \omega x \omega^{-1}, \qquad x\in H, \\
		\label{eq:omega-traza}
		\tr_V(\theta \omega) &= \tr_V(\theta \omega^{-1}),  \qquad  \theta \in \End_H(V),\ V\in \Rep H.
	\end{align}
	When $\omega\in G(H)$ satisfies~\eqref{eq:omega-square-antipode}, it is called a \emph{pivot}. If, in addition, $\omega$ fulfills~\eqref{eq:omega-traza}, 
	then it is a \emph{spherical element}. Observe that this is automatic if $\omega^2=1$. See~\cite{A+} and references therein for notation and further details. 
	
	In this setting, it is possible to define a \emph{quantum trace} for each $V\in \Rep H$ and $\theta \in \End_H(V)$, via
	$\qtr(\theta)=\tr(\theta\omega)$, that is, the trace of the map $V\to V$, $v\mapsto \theta(\omega\cdot v)$.
	In turn, this gives rise to a \emph{quantum dimension}:
	\begin{align*}
		\qdim\colon &\Rep H\to \k, & \qdim V= \qtr(\id_V)=\tr_V(\omega).
	\end{align*}
	A quotient category $\uRep H$ is defined, with objects $\{[X] : X\in \Rep H\}$ and 
	morphisms 
	$\underline{\Hom}([X], [Y]) \coloneqq  \Hom(X, Y)/J(X,Y)$, $X,Y \in \Rep H$, where
	\begin{align*}
		J(X,Y) &= \{f \in  \Hom(X, Y): \qtr(fg) = 0, \forall g\in  \Hom(Y, X) \}.
	\end{align*}
	The category $\uRep H$ is semisimple. The irreducible objects are the classes $[M]$ of the indecomposable $H$-modules $M$ with $\qdim M\neq0$;  see~\cite{BaW-adv}.
	
	\subsection{$H_\lambda$ is spherical}
	From now on, we assume that $N$ is odd. Recall the Hopf algebra structure of $H_\lambda$ in~\eqref{eqn:comultiplication}.
	
	\begin{lemma}
		$H_\lambda$ is a spherical Hopf algebra, with pivot $\omega=g_1^N$.
	\end{lemma}
	\pf
	One can verify directly that $S^2(x)=g_1^Nxg_1^{-N}$, for any $x\in H_\lambda$; it holds trivially on $\Gamma$, and it is enough to check it on the generators $a_1$ and $a_2$. The pivot is a spherical element because it is an involution.
	\epf
	
	\begin{example}\label{exa:qdim-simple}
		Let $L$ be an irreducible $H_\lambda$-module. If $\dim L=1$, then $L=L_1^\chi$, with $\chi=(\zeta^i,\xi^j)\in\O_1$. We have $[g_1^N]=(\zeta^{Ni})=\pm1$ and thus $\qdim L=\pm1$.
		On the other hand, if $\dim L>1$, then $\qdim L=0$. Indeed, if $\dim L=2$, then $L=L_{2,h}^\chi$ or $L_{2,v}^\chi$, with $\chi\in\O_2$, and $[g_1^N]=\diag(\zeta^{Ni},-\zeta^{Ni})$; if $\dim L=4$, then $[g_1^N]=\diag(\zeta^{Ni},-\zeta^{Ni},\zeta^{Ni},-\zeta^{Ni})$.
	\end{example}
	
	If we combine Lemma \ref{lem:indecomp} with Example \ref{exa:qdim-simple} above, then we get:
	\begin{corollary}\label{cor:dim-O2}
		If $M=M^\chi$ is indecomposable and $\chi\in\O_2\sqcup\O_4$, then $\qdim M=0$.\qed
	\end{corollary}
	
	In turn, Proposition \ref{pro:extensions} and Lemma \ref{lem:indecomp4} give:
	\begin{corollary}
		If $M$ is an indecomposable $H_\lambda$-module of dimension 2 or 4, then $\qdim M=0$.\qed
	\end{corollary}
	
	\begin{remark}
		Corollary \ref{cor:dim-O2} implies  that $\qdim Q_{n,h}^\chi=\qdim Q_{n,v}^\chi=0$ for the modules in Definition \ref{def:indecn-2}.
		In turn, for the indecomposable modules $Q_{n}^\chi$ from Definition \ref{def:indecn-1}, we get $\qdim Q_n^\chi=\zeta^{Ni}=\pm1$ when $n$ is odd, where $\chi=(\zeta^i,\xi^j)\in\O_1$, and $\qdim Q_n^\chi=0$ when $n$ is even.
	\end{remark}
	
	We expect the answer to
	the following question to be affirmative.
	\begin{question}\label{question:zero}
		Let $M$ be an indecomposable $H_\lambda$-module of even dimension. 
		Is $\qdim M=0$?
	\end{question}

	\subsection{Lower fusion rules}
	
	We compute tensor products between irreducible objects in $\uRep H_\lambda$ coming from  $H_\lambda$-modules of dimension 3, cf.~\S\ref{sec:small-dim}.
	
	\begin{proposition}
		Let $\chi,\phi\in\O_1$.
		\begin{enumerate}[leftmargin=*]
			\item The following are indecomposable and  pairwise non-isomorphic:
			\begin{align*}
				M_{3,1}^{\chi}&\ot M_{3,1}^{\phi},&& M_{3,1}^{\chi}\ot M_{3,2}^{\phi}, &&M_{3,1}^{\chi}\ot M_{3,3}^{\phi}, && M_{3,2}^{\chi}\ot M_{3,2}^{\phi},\\ M_{3,2}^{\chi}&\ot M_{3,4}^{\phi}, && M_{3,3}^{\chi}\ot M_{3,3}^{\phi}, &&M_{3,3}^{\chi}\ot M_{3,4}^{\phi}, && M_{3,4}^{\chi}\ot M_{3,4}^{\phi}. 
			\end{align*}
			Also, $M^{\chi}_{3,i}\ot M^{\phi}_{3,j}\simeq M^{\chi}_{3,j}\ot M^{\phi}_{3,i}$ for $(i,j)\in\{(2,1),(3,1),(3,4),(4,2)\}$.
			\item We have $M_{3,1}^{\chi}\ot M_{3,4}^{\phi}\simeq P^{\chi\phi}\oplus L_1^{-\overline{\chi\phi}}\simeq M_{3,4}^{\chi}\ot M_{3,1}^{\phi}$ and
			\begin{align*}
				M_{3,2}^{\chi}\ot M_{3,3}^{\phi}&\simeq C_{2,-1}^{-\overline{\chi\phi}}\oplus C_{2,1}^{\chi\phi}\oplus L_{1}^{\overline{\chi\phi}}, &
				M_{3,3}^{\chi}\ot M_{3,2}^{\phi}&\simeq C_{2,1}^{-\overline{\chi\phi}}\oplus C_{2,-1}^{\chi\phi}\oplus L_{1}^{\overline{\chi\phi}}.
			\end{align*}
		\end{enumerate}
		In the category $\uRep H_\lambda$, the tensor products decompose as: $[M_{3,1}^{\chi}]\ot [M_{3,4}^{\phi}]\simeq [M_{3,4}^{\chi}]\ot [M_{3,1}^{\phi}]\simeq [L_1^{-\overline{\chi\phi}}]$ and
		$[M_{3,2}^{\chi}]\ot [M_{3,3}^{\phi}]\simeq [M_{3,3}^{\chi}]\ot [M_{3,2}^{\phi}]\simeq [L_{1}^{\overline{\chi\phi}}]$ while the other tensor products give rise to new simple objects.
	\end{proposition}
	\begin{remark}\label{rem:notiso}
		Observe that $M_{3,2}^{\chi}\ot M_{3,3}^{\phi}\not \simeq M_{3,3}^{\phi}\ot M_{3,2}^{\chi}$. 
		In particular, $H_\lambda$ is not quasi-triangular. 
	\end{remark}
	\pf
	$(1)$ We claim that for each $i\in\I_4$ there is a basis $\{x_\ell:\ell\in\I_9\}$ so that $M_{3,i}^{\chi}\ot M_{3,i}^{\phi}$ is described by the following diagrams, respectively:
	\begin{align*}		
		\resizebox{3cm}{.15cm}
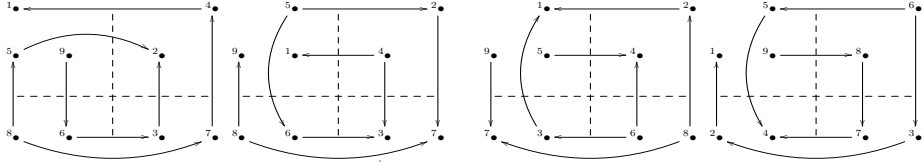
{
			\xymatrix{
				^{1}\bullet&&\ar@{--}[ddd]&&^{4}\bullet\ar[llll]\\
				^{5}\bullet\ar@/^1.5pc/[rrr]&^{9}\bullet\ar[dd]&&^{2}\bullet&\\
				\ar@{--}[rrrr]&&&&\\
				^{8}\bullet\ar@/_1.5pc/[rrrr]\ar[uu]&^{6}\bullet\ar[rr]&&^{3}\bullet\ar[uu]&^{7}\bullet\ar[uuu]
			}
		}
		&
		\resizebox{3cm}{.15cm}{
			\xymatrix{
				&^{5}\bullet\ar[rrr]\ar@/_1.5pc/[ddd]&\ar@{--}[ddd]&&^{2}\bullet\ar[ddd]\\
				^{9}\bullet&^{1}\bullet&&^{4}\bullet\ar[ll]\ar[dd]&\\
				\ar@{--}[rrrr]&&&&\\
				^{8}\bullet\ar@/_1.5pc/[rrrr]\ar[uu]&^{6}\bullet\ar[rr]&&^{3}\bullet&^{7}\bullet
			}
		}
		&
		\resizebox{3cm}{.15cm}{
			\xymatrix{
				&^{1}\bullet&\ar@{--}[ddd]&&^{2}\bullet\ar[lll]\\
				^{9}\bullet\ar[dd]&^{5}\bullet\ar[rr]&&^{4}\bullet&\\
				\ar@{--}[rrrr]&&&&\\
				^{7}\bullet&^{3}\bullet\ar@/^1.5pc/[uuu]&&^{6}\bullet\ar[uu]\ar[ll]&^{8}\bullet\ar[uuu]\ar@/^1.5pc/[llll]\\
			}
		}
		&
		\resizebox{3cm}{.15cm}{
			\xymatrix{
				&^{5}\bullet\ar@/_1.5pc/[ddd]&\ar@{--}[ddd]&&^{6}\bullet\ar[ddd]\ar[lll]\\
				^{1}\bullet&^{9}\bullet\ar[rr]&&^{8}\bullet\ar[dd]&\\
				\ar@{--}[rrrr]&&&&\\
				^{2}\bullet\ar[uu]&^{4}\bullet&&^{7}\bullet\ar[ll]&^{3}\bullet\ar@/^1.5pc/[llll]\\
			}
		}
	\end{align*}
	The tensor products $M_{3,i}^\chi\ot M_{3,j}^\phi$, $i\neq j$, can also be described via diagrams, which lead to the isomorphisms from the statement. Explicitly, the modules $M_{3,i}^\chi\ot M_{3,j}^\phi$, $(i,j)=(1,2), (1,3), (2,4)$ and $(3,4)$ are, respectively:
	\begin{align*}
		\resizebox{3cm}{.15cm}
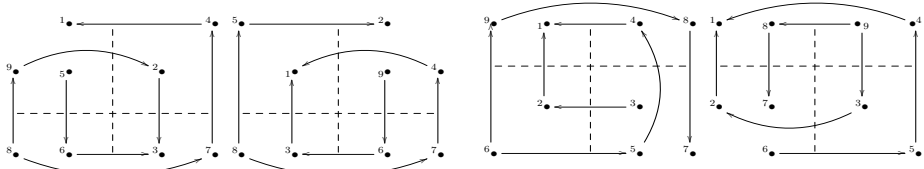
{
			\xymatrix{
				&^{1}\bullet&\ar@{--}[ddd]&&^{4}\bullet\ar[lll]\\
				^{9}\bullet\ar@/^1.5pc/[rrr]&_{5}\bullet\ar[dd]&&^{2}\bullet\ar[dd]&\\
				\ar@{--}[rrrr]&&&&\\
				^{8}\bullet\ar[uu]\ar@/_1.5pc/[rrrr]&^{6}\bullet\ar[rr]&&^{3}\bullet&^{7}\bullet\ar[uuu]
			}
		}
		&
		\resizebox{3cm}{.15cm}{
			\xymatrix{
				^{5}\bullet\ar[rrr]&&\ar@{--}[ddd]&^{2}\bullet&\\
				&_{1}\bullet&&_{9}\bullet\ar[dd]&^{4}\bullet\ar@/_1.5pc/[lll]\\
				\ar@{--}[rrrr]&&&&\\
				^{8}\bullet\ar[uuu]\ar@/_1.5pc/[rrrr]&^{3}\bullet\ar[uu]&&^{6}\bullet\ar[ll]&^{7}\bullet\ar[uu]
			}
		}
		&
		\resizebox{3cm}{.15cm}{
			\xymatrix{
				^{9}\bullet\ar@/^1.5pc/[rrrr]&_{1}\bullet&\ar@{--}[ddd]&^{4}\bullet\ar[ll]&^{8}\bullet\ar[ddd]\\
				\ar@{--}[rrrr]&&&&\\
				&^{2}\ar[uu]\bullet&&^{3}\bullet\ar[ll]&\\
				^{6}\bullet\ar[uuu]\ar[rrr]&&&^{5}\bullet\ar@/_1.5pc/[uuu]&^{7}\bullet
			}
		}
		&
		\resizebox{3cm}{.15cm}{
			\xymatrix{
				^{1}\bullet&_{8}\bullet\ar[dd]&\ar@{--}[ddd]&\bullet\,_{9}\ar[ll]\ar[dd]&\bullet^{4}\ar@/_1.5pc/[llll]\\
				\ar@{--}[rrrr]&&&&\\
				^{2}\bullet\ar[uu]&^{7}\bullet&&^{3}\bullet\ar@/^1.5pc/[lll]&\\
				&^{6}\bullet\ar[rrr]&&&^{5}\bullet\ar[uuu]
			}
		}
	\end{align*}
	We work out the case $M=M^\chi_{3,2}\ot M^\phi_{3,4}$; the other cases follow similarly. 
	Recall from Definition \ref{def:indec3} that $M^\chi_{3,2}$ is the $H_\lambda$-module with basis $\{u,v,w\}$, determined by the facts that $u\in M^\chi_{3,2}[\chi]$ and the only non-trivial actions are $a_1\cdot v=u$, $a_2\cdot v=w$. Similarly, $M^\phi_{3,4}$ has a basis $\{u',v',w'\}$ such that $u'\in M^\phi_{3,4}[\phi]$ and $a_1\cdot w'=v'$, $a_2\cdot v'=u'$. We consider the natural basis $\{x_\ell:\ell\in\I_9\}$ of $M_{3,2}^{\chi}\ot M_{3, 4}^{\phi}$, with
	\begin{align*}
		x_1 &=u\ot u', && x_2 =u\ot v' ,&&  x_3 =u\ot w' ,&&  x_4 =v\ot u' ,&&  x_5 =v\ot v' ,\\
		x_6 &=v\ot w' , && x_7 =w\ot u' , && x_8 =w\ot v' ,&&  x_9 =w\ot w' .
	\end{align*}
	We make the following change of basis:
	\begin{align*}
		\tilde{x}_1&=-\chi_1\chi_2 x_1 ,&&\tilde{x}_2= -\chi_1 x_2, &&\tilde{x}_3= -x_3,&&\tilde{x}_4= -\chi_1 (x_8 + \chi_2 x_4), \\ \tilde{x}_5&= x_3 - \chi_1 x_5 ,&&
		\tilde{x}_6= x_6 , &&\tilde{x}_7= -\chi_1\chi_2 x_7 ,&&\tilde{x}_8=  \chi_1 x_8,\quad	\tilde{x}_9= x_9,
	\end{align*}
	which leads to the module described above. 
	
	To check indecomposability, we follow the ideas of \S\ref{sec:any-dim}. 
	Observe that the simple submodules of $M$ are $L_1^{-\overline{\chi\phi}}\simeq \lg x_7\rg\subset M$ and $L_1^{\chi\phi}\simeq \lg x_1\rg\subset M$. 
	In particular, if we write $M=N_1\oplus N_2$ and assume that $x_7\in N_1$, then this forces $x_8\in N_1$, so $x_9\in N_1$. 
	If $N_2\neq\{0\}$, then $x_1\in N_2$ and thus $x_4\in N_2$. Since $x_9\in N_1$, we necessarily have $x_6\in N_1$ (and $x_2\in N_2$). This forces $x_5\in N_1$  and consequently $x_4\in N_1$, a contradiction. Hence $M=N_1$.
	
	$(2)$ The module $M_{3,1}^\chi$ is generated by an element $w\in M_{3,1}^\chi[-\overline{\chi}]$ and has basis $\{w,v\coloneqq a_2\cdot w,u\coloneqq a_1a_2\cdot w\}$. Similarly, $M_{3,4}^\phi$ is generated by $w'\in M_{3,4}^\phi[-\overline{\phi}]$ and $\{w',v'\coloneqq a_1\cdot w',u'\coloneqq a_2a_1\cdot w'\}$ is a basis.
	We check that $a_1a_2a_1a_2\cdot (w\ot w')\neq 0$; hence it determines a direct summand $P^{\chi\phi}$ and a simple complement 
	$\lg u\ot w'-\chi_1 v\ot v'+\chi_1\chi_2 w\ot u' \rg\simeq L_1^{-\overline{\chi\phi}}$.
	
	The remaining decompositions follow as in $(1)$.
	\epf

\end{document}